\documentclass[11pt,a4paper]{article}
\usepackage{amsfonts,amsmath,amssymb,accents,bm,amsthm}
\usepackage{verbatim}

\usepackage{hyperref}
\usepackage[normalem]{ulem}

\usepackage{diagbox}

\usepackage{graphicx}
\usepackage{xcolor}

\newcommand{\dis}{\displaystyle}

\newcommand{\noi}{\noindent}
\newcommand{\halmos}{\rule{1ex}{1.4ex}}
\newcommand{\QED}{\nopagebreak{\hspace*{\fill}$\halmos$\medskip}}

\newcommand{\quand}{\quad\mbox{and}\quad}

\newtheoremstyle{mythm}
  {}
  {}
  {\itshape}
  {}
  {\bfseries}
  {}
  {.5em}
  {#1 #2 \thmnote{(#3)}}

\theoremstyle{mythm}
\newtheorem{theorem}{Theorem}
\newtheorem{proposition}[theorem]{Proposition}
\newtheorem{lemma}[theorem]{Lemma}

\newtheorem{exercise}[theorem]{Exercise}
\newtheorem{corollary}[theorem]{Corollary}
\newtheorem{conjecture}[theorem]{Conjecture}

\newtheorem{counterex}[theorem]{Counterexample}

\newcommand{\bt}{\begin{theorem}}
\newcommand{\et}{\end{theorem}}
\newcommand{\bl}{\begin{lemma}}
\newcommand{\el}{\end{lemma}}
\newcommand{\bp}{\begin{proposition}}
\newcommand{\ep}{\end{proposition}}
\newcommand{\bcor}{\begin{corollary}}
\newcommand{\ecor}{\end{corollary}}
\newcommand{\br}{\begin{remark}\rm}
\newcommand{\er}{\end{remark}}
\newcommand{\bcon}{\begin{conjecture}}
\newcommand{\econ}{\end{conjecture}}
\newcommand{\bex}{\begin{exercise}}
\newcommand{\eex}{\end{exercise}}
\newcommand{\bcou}{\begin{counterex}}
\newcommand{\ecou}{\end{counterex}}

\theoremstyle{definition}

\newtheorem{remark}[theorem]{Remark}

\newenvironment{Proof}[1][]{\noi\textbf{Proof #1}}{\QED}
\newcommand{\bpro}{\begin{Proof}}
\newcommand{\epro}{\end{Proof}}

\newcommand{\be}{\begin{equation}}
\newcommand{\ee}{\end{equation}}
\newcommand{\ba}{\begin{array}}
\newcommand{\ea}{\end{array}}
\newcommand{\bac}{\begin{array}{r@{\,}c@{\,}l}}
\newcommand{\bc}{\be\begin{array}{r@{\,}c@{\,}l}}
\newcommand{\ec}{\end{array}\ee}

\newcommand{\De}{\Delta}

\newcommand{\La}{\Lambda}

\newcommand{\Fi}{{\cal F}}

\newcommand{\Hi}{{\cal H}}

\newcommand{\Li}{{\cal L}}
\newcommand{\Mi}{{\cal M}}

\newcommand{\Ri}{{\cal R}}
\newcommand{\Si}{{\cal S}}

\newcommand{\E}{{\mathbb E}}
\newcommand{\F}{{\mathbb F}}

\renewcommand{\P}{{\mathbb P}}

\newcommand{\R}{{\mathbb R}}

\newcommand{\Rb}{{\mathbf R}}
\newcommand{\Sb}{{\mathbf S}}

\newcommand{\Xb}{{\mathbf X}}
\newcommand{\Yb}{{\mathbf Y}}

\newcommand{\sub}{\subset}
\newcommand{\beh}{\backslash}

\newcommand{\ti}{\tilde}

\newcommand{\dgg}{\dagger}

\newcommand{\un}{\underline}

\newcommand{\eb}{\mathbf{e}}

\newcommand{\xb}{\mathbf{x}}
\newcommand{\yb}{\mathbf{y}}

\newcommand{\fb}{\mathbf{f}}
\newcommand{\gb}{\mathbf{g}}

\newcommand{\psib}{{\boldsymbol\psi}}

\begin{document}

\makeatletter\@addtoreset{equation}{section}
\makeatother\def\theequation{\thesection.\arabic{equation}}

\renewcommand{\labelenumi}{{\rm (\roman{enumi})}}
\renewcommand{\theenumi}{\roman{enumi}}

\title{\vspace*{-2cm}Commutative monoid duality}
\author{Jan Niklas Latz\footnote{The Czech Academy of Sciences,
  Institute of Information Theory and Automation,
  Pod vod\'arenskou v\v e\v z\' i~4,
  18200 Praha 8.
  Czech Republic.}\:\,\footnote{latz@utia.cas.cz.}
\and
Jan~M.~Swart${}^\ast$\footnote{swart@utia.cas.cz, ORCID: 0000-0001-8614-4053.}}

\date{\today}

\maketitle

\begin{abstract}\noi
We introduce two partially overlapping classes of pathwise dualities between interacting particle systems that are based on commutative monoids (semigroups with a neutral element) and semirings, respectively. For interacting particle systems whose local state space has two elements, this approach yields a unified treatment of the well-known additive and cancellative dualities. For local state spaces with three or more elements, we discover several new dualities.
\end{abstract}
\vspace{.5cm}

\noi
{\it MSC 2020.} Primary: 82C22; Secondary: 16Y60, 20M32. \\
{\it Keywords:} interacting particle system, duality, monoid, semiring. \\
{\it Acknowledgements:} Work supported by grant 20-08468S of the Czech Science Foundation (GA\v{C}R). \\
%
%
{\it Data Availability Statement.} Data sharing not applicable to this article as no datasets were generated or analysed during the current study.

\newpage

{\setlength{\parskip}{-2pt}\tableofcontents}

\newpage

\section{Introduction}\label{S:intro}

\subsection{Aim of the paper}

The use of duality in the study of Markov processes in general, and of interacting particle systems in particular, has a long history see \cite[Section~II.3]{Lig85}. While in the past, many useful dualities have been discovered by ad hoc methods, more recently several authors have attempted a more systematic search. There are two approaches: the pathwise approach propagated in, e.g., \cite{CS85,JK14,SS18}, and the algebraic approach of \cite{LS95,Sud00,GKRV09,CGGR15,SSF18}. In the present paper, we use the pathwise approach to treat a number of known dualities in a unified framework and discover new dualities.

Let $\Rb,\Sb$, and $T$ be sets and let $m:\Sb\to\Sb$, $n:\Rb\to\Rb$, and $\psib:\Sb\times\Rb\to T$ be functions. By definition, we say that the map $m$ is \emph{dual} to $n$ with respect to the \emph{duality function} $\psib$ if
\be\label{mapdual}
\psib\big(m(x),y\big)=\psib\big(x,n(y)\big)\qquad(x\in\Sb,\ y\in\Rb).
\ee
We will especially be interested in the case that $\Sb$ and $\Rb$ are product spaces of the form $\Sb=S^\La$ and $\Rb=R^\La$ where $S,R$, and $\La$ are finite sets. We will be interested in finding triples $(S,R,T)$ for which the following is true:
\begin{quote}
For each finite set $\La$, there exists a function $\psib:S^\La\times R^\La\to T$, a class $\Si$ of functions $m:S^\La\to S^\La$, and a class $\Ri$ of functions $n:R^\La\to R^\La$, so that each $m\in\Si$ has a unique dual map $n\in\Ri$ with respect to $\psib$, and vice versa.
\end{quote}
Recall that a monoid is a semigroup that contains a neutral element. For our first result (the combination of Propositions \ref{P:mondu} and \ref{P:prodrefl} below), we will assume that $S,R$, and $T$ are commutative monoids. In this case $R=S^\La$ and $S=R^\La$ also naturally have the structure of a monoid and $\Si$ and $\Ri$ will be the classes of monoid homomorphisms from $S$ to $S$ and from $R$ to $R$, respectively. For our second result (Proposition~\ref{P:lindu} below), we will assume that $S=R=T$ is a semiring. In this case $S^\La$ has the structure of a left module over $S$ and also the structure of a right module over $S$. Now $\Si$ and $\Ri$ will be the classes of maps $m:S^\La\to S^\La$ that preserve the structure of $S^\La$ as a left or right module, respectively.

In Section~\ref{S:numeric}, we will explicitly find all triples $(S,R,T)$ with $S,R$, and $T$ having at most four elements that satisfy the conditions of our main results (Propositions \ref{P:mondu} and \ref{P:lindu}). There is a considerable overlap, in the sense that several of the triples $(S,R,T)$ that satisfy the conditions of Proposition~\ref{P:mondu} also satisfy the conditions of Proposition~\ref{P:lindu}, and vice versa. While our main results are algebraic in nature, our motivation comes from probability theory. We spend the remainder of this section explaining the motivation of our work and then turn to the purely algebraic questions.

\subsection{Pathwise duality of Markov processes}\label{S:path}

Our main motivation comes from the theory of interacting particle systems, as we now explain. Let $\Sb$ be a finite set (typically of the form $\Sb=S^\La$), let $\Mi$ be a finite set whose elements are maps $m:\Sb\to\Sb$, and let $(r_m)_{m\in\Mi}$ be nonnegative real constants. One is frequently interested in continuous-time Markov processes $(X_t)_{t\geq 0}$ that evolve according to the following informal description:
\begin{itemize}
\item At the times of a Poisson process with intensity $r_m$, the previous state $x$ of the process is replaced by the new state $m(x)$.
\end{itemize}

More precisely, such a process can be constructed as follows. Let $\mu$ be the measure on $\Mi$ defined by $\mu(\{m\}):=r_m$, let $\ell$ denote the Lebesgue measure on $\R$, let $\mu\otimes\ell$ denote the product measure on $\Mi\times\R$, and let $\Pi$ be a Poisson subset of $\Mi\times\R$ with intensity measure $\mu\otimes\ell$. For each $s,u\in\R$ with $s\leq u$, define
\be\ba{ll}
\dis\Pi^+_{s,u}:=\big\{(m,t)\in\Pi:s<t\leq u\big\},\quad
\dis\Pi^-_{s,u}:=\big\{(m,t)\in\Pi:s\leq t<u\big\},
\ec
and define random maps by
\be
\Xb^\pm_{s,u}:=m_n\circ\cdots\circ m_1
\quad\mbox{where}\quad
\Pi^\pm_{s,u}=\big\{(m_1,t_t),\ldots,(m_n,t_n)\big\}
\quad\mbox{with}\quad
t_1<\cdots<t_n.
\ee
The collection of random maps $(\Xb^\pm_{s,u})_{s\leq u}$ is called a \emph{stochastic flow}. It is easy to see that $\Xb^\pm_{s,s}$ is the identity map and that $\Xb^\pm_{t,u}\circ\Xb^\pm_{s,t}=\Xb^\pm_{s,u}$ for all $s\leq t\leq u$. 

Let $X_0$ be an $\Sb$-valued random variable, independent of the Poisson set $\Pi$. Then one can prove that for each $s\in\R$, setting
\be
X^\pm_t:=\Xb^\pm_{s,s+t}(X_0)\qquad(t\geq 0)
\ee
defines Markov processes $(X^-_t)_{t\geq 0}$ and $(X^+_t)_{t\geq 0}$ that fit our informal description above. The process $(X^-_t)_{t\geq 0}$ has left-continuous sample paths while $(X^+_t)_{t\geq 0}$ has right-continuous sample paths. Since at deterministic times, $X^-_t=X^+_t$ almost surely, in practice it does not matter too much which version of the process we use.

Now assume that $\Rb$ and $T$ are finite sets, that $\psib:\Sb\times\Rb\to T$ is a function, and that each map $m\in\Mi$ has a unique dual map $\hat m$ with respect to $\psib$. Let $\hat\Mi:=\{\hat m:m\in\Mi\}$. Then setting 
\be
\hat\Pi:=\big\{(\hat m,-t):(m,t)\in\Pi\big\}
\ee
defines a Poisson subset of $\hat\Mi\times\R$. We can use $\hat\Pi$ to construct stochastic flows $(\Yb^-_{s,u})_{s\leq u}$ and $(\Yb^+_{s,u})_{s\leq u}$ precisely in the same way as we did for the Poisson set $\Pi$. Then it is easy to see that
\be\label{pathdual}
\psib\big(\Xb^\pm_{s,u}(x),y\big)=\psib\big(x,\Yb^\mp_{-u,-s}(y)\big)
\qquad(s\leq u,\ x\in\Sb,\ y\in\Rb),
\ee
i.e., the map $\Xb^+_{s,u}$ is dual to $\Yb^-_{-u,-s}$ and likewise $\Xb^-_{s,u}$ is dual to $\Yb^+_{-u,-s}$. In the theory of Markov processes, a duality relation between stochastic flows of the form (\ref{pathdual}) is called a \emph{pathwise duality}.

\subsection{Additive and cancellative duality}

We will especially be interested in \emph{interacting particle systems}, which are Markov processes with a state space of the form $\Sb=S^\La$, where $S$ is a finite set, called the \emph{local state space}, and $\La$ is any finite or countably infinite set that is usually called the \emph{lattice} (not to be confused with the order theoretic lattices that we will discuss later). Elements of $S^\La$ are functions $\xb:\La\to S$. For technical simplicity, we will only discuss finite $\La$.

Two forms of duality, called \emph{additive} and \emph{cancellative} duality, have found widespread applications in the theory of interacting particle systems \cite{Gri79,Lig85}. To explain these, let $\La$ be a finite set and let $\Sb=\Rb:=\{0,1\}^\La$ and $T:=\{0,1\}$. Let $\psib_{\rm add}:\Sb\times\Sb\to T$ and $\psib_{\rm canc}:\Sb\times\Sb\to T$ be defined by:
\be\label{adccan}
\psib_{\rm add}(\xb,\yb):=\bigvee_{i\in\La}\xb(i)\yb(i)
\quand
\psib_{\rm canc}(\xb,\yb):=\sum_{i\in\La}\xb(i)\yb(i)\ {\rm mod}(2).
\ee
One can prove that a map $m:\Sb\to\Sb$ has a dual with respect to $\psib_{\rm add}$ if and only if it is \emph{additive}, which means that
\be
m(\un 0)=0\quand m(\xb\vee\yb)=m(\xb)\vee m(\yb)\qquad(\xb,\yb\in\Sb),
\ee
where $\un 0(i):=0$ $(i\in\La)$ denotes the function that is identically zero. Similarly, a map $m:\Sb\to\Sb$ has a dual with respect to $\psib_{\rm canc}$ if and only if it is \emph{cancellative}, which means that
\be
m(\un 0)=0\quand m\big(\xb+\yb\ {\rm mod}(2)\big)=m(\xb)+m(\yb)\ {\rm mod}(2)\qquad(\xb,\yb\in\Sb).
\ee
The duals of additive or cancellative maps, if they exist, are unique and such dual maps are also additive or cancellative, respectively. A Markov process is called \emph{additive} or \emph{cancellative} if it can be constructed using only maps of the appropriate type. Some of the most studied interacting particle systems are additive, including the voter model, the contact process, and the exclusion process \cite{Lig99}, and duality is one of the most important tools in their study. Cancellative duality has succesfully been applied in the study of various nonlinear voter models \cite{CD91,Han99,SS08} and annihilating branching processes \cite{BDD91}. We are motivated by the wish to find generalisations of the duality functions in (\ref{adccan}) to local state spaces $S$ with three or more elements.

\subsection{A new form of duality}

As an appetizer for the remainder of the paper, we highlight one particular duality that we have found as a consequence of our results. Let $S:=\{0,1,2\}$ be equipped with the binary operation $\oplus$ that is defined by the following addition table:

\begin{center}
\begin{tabular}{c|ccc}
$\oplus$&0&1&2\\
\hline
0&0&1&2\\
1&1&2&1\\
2&2&1&2\\
\end{tabular}
\end{center}

Let $R=T:=\{-1,0,1\}$, equipped with the usual product. Then one can check that $S$ and $R$ are commutative monoids. Indeed, in Subsection~\ref{S:examp} below, we list all commutative monoids with at most three elements. In the notation used there, $S=M_6$ and $R\cong M_5$. We now fix a finite set $\La$ and for $x\in S$, we let $\un x\in S^\La$ denote the function that is constantly $x$, i.e., $\un x(i):=x$ $(i\in\La)$. For $y\in R$, we define $\un y\in R^\La$ similarly. For $\xb,\yb\in S^\La$, we define $\xb\oplus\yb$ in a pointwise way, i.e., $(\xb\oplus\yb)(i):=\xb(i)\oplus\yb(i)$ $(i\in\La)$. For $\xb,\yb\in R^\La$, we define the pointwise product $\xb\cdot\yb$ similarly. We define $\Si$ to be the set of all functions $m:S^\La\to S^\La$ such that
\be
m(\un 0)=0\quand m(\xb\oplus\yb)=m(\xb)\oplus m(\yb)\qquad(\xb,\yb\in S^\La).
\ee
Similarly, we let $\Ri$ denote the set of all functions $n:R^\La\to R^\La$ such that
\be
n(\un 1)=1\quand n(\xb\cdot\yb)=n(\xb)\cdot n(\yb)\qquad(\xb,\yb\in R^\La).
\ee
We define $\psi:S\times R\to R$ by
\be
\left(\ba{ccc}
\psi(0,-1)&\psi(0,0)&\psi(0,1)\\
\psi(1,-1)&\psi(1,0)&\psi(1,1)\\
\psi(2,-1)&\psi(2,0)&\psi(2,1)\ea\right)
:=
\left(\ba{ccc}
1&1&1\\
-1&1&1\\
0&0&1\ea\right)
\ee
which corresponds to the function $\psi_5$ from Subsection~\ref{S:moex}, and we define $\psib:S^\La\times R^\La\to R^\La$ by
\be\label{new}
\psib(\xb,\yb):=\prod_{i\in\La}\psi\big(\xb(i),\yb(i)\big)\qquad(\xb\in S^\La,\ \yb\in R^\La).
\ee
Then as an immediate consequence of Propositions \ref{P:mondu} and \ref{P:prodrefl} below, we obtain the following result.

\bp[A new duality]
Each\label{P:new} map $m\in\Si$ has a unique dual map $\hat m$ with respect to the duality function $\psib$ defined in (\ref{new}), and this dual map satisfies $\hat m\in\Ri$. Conversely, for each $n\in\Ri$, there exists a unique $m\in\Si$ such that $n$ is the dual of $m$ with respect to $\psib$.
\ep

Let $\Mi$ be a subset of $\Si$, and let $(r_m)_{m\in\Mi}$ be nonnegative rates. Let $(X_t)_{t\geq 0}$ be an interacting particle system that is constructed by applying each map $m\in\Mi$ at the times of a Poisson process with intensity $r_m$. By the general principles explained in Subsection~\ref{S:path}, such an interacting particle system is pathwise dual to an interacting particle system $(Y_t)_{t\geq 0}$ that is constructed by applying each dual map $\hat m$ at the times of Poisson process with intensity $r_m$. Letting $(X^x_t)_{t\geq 0}$ and $(Y^y_t)_{t\geq 0}$ denote the processes started in the initial states $X^\xb_0=\xb$ and $Y^\yb_0=\yb$, it is easy to see that
\be\label{dual}
\E\big[\psib(X^\xb_t,\yb)\big]=\E\big[\psib(\xb,Y^\yb_t)\big]
\qquad\big(\xb\in S^\La,\ \yb\in R^\La,\ t\geq 0).
\ee
Indeed, this follows by setting $X^\xb_t:=\Xb^+_{0,t}(x)$ and $Y^\yb_t:=\Yb^-_{-t,0}(y)$ and taking expectations in (\ref{pathdual}), using the fact that $R$ is naturally embedded in $\R$.

\subsection{Open problems}

For all duality functions $\psib$ considered in this paper, it will be true that knowing $\psib(\xb,\yb)$ for all $\yb\in R^\La$ uniquely determines $\xb\in S^\La$. As a consequence, if $(\Xb^\pm)_{s\leq u}$ and $(\Yb^\pm)_{s\leq u}$ are dual stochastic flows as in Subsection~\ref{S:path}, then the law of $\Xb^\pm_{0,t}(x)$ is uniquely determined by all probabilities of the form
\be\label{fdd}
\P\big[\psib\big(\Xb^\pm_{0,t}(\xb),\yb_1\big)=z_1,\ldots,\psib\big(\Xb^\pm_{0,t}(\xb),\yb_n\big)=z_n\big],
\ee
with $\yb_1,\ldots,\yb_n\in R^\La$ and $z_1,\ldots,z_n\in T$. By the pathwise duality relation (\ref{pathdual}), the probability in (\ref{fdd}) equals
\be
\P\big[\psib\big(\xb,\Yb^\mp_{-t,0}(\yb_1)\big)=z_1,\ldots,\psib\big(\xb,\Yb^\mp_{-t,0}(\yb_n)\big)=z_n\big].
\ee
For the duality highlighted in Proposition~\ref{P:new}, the situation turns out to be considerably better. In fact, in this example, one can prove that the law of $\Xb^\pm_{0,t}(x)$ is uniquely determined by all expectations of the form
\be\label{Epsi}
\E\big[\psib\big(\Xb^\pm_{0,t}(\xb),\yb\big)\big]
\ee
with $\yb\in R^\La$. In general, it is not hard to see that each finite monoid $T$ can be represented in a real algebra. One can then view $\psib$ as a function taking values in this real algebra and define expectations as in (\ref{Epsi}). However, in this generality it is not true for all dualities that we will find in the sections to come that the law of $\Xb^\pm_{0,t}(x)$ is uniquely determined by all expectations of the form (\ref{Epsi}). Therefore, we pose as an open problem to classify all dualities for which distributional uniqueness holds in this stronger form. A more vaguely formulated problem is to determine more generally the ``minimal'' information one needs about probabilities of the form (\ref{fdd}) to determine the law of $\Xb^\pm_{0,t}(x)$ uniquely.

Another vaguely formulated open problem concerns further generalisations of our results. Our main results, Propositions \ref{P:mondu} and \ref{P:lindu}, are in many ways similar, which leads one to suspect it may be possible to combine them into one even more general (but presumably even more abstract) result. At present, we do not know how this should be done.

Finally, we note that Lloyd and Sudbury \cite{LS95,Sud00} have studied general duality functions that can be written as a product over the set $\La$ as in (\ref{new}). The work in \cite{LS95,Sud00} is restricted to local state spaces with two elements. They have found useful dualities of the form (\ref{dual}) that do not always come from pathwise dualities of the form (\ref{pathdual}) and that in some way interpolate between additive and cancellative duality (see also \cite[Section~2.7]{Swa13}). Our present work was motivated by the wish to generalise their work to state spaces with three and more elements. However, we still do not know if there is an elegant way to do this.

\subsection{Outline}

The outline of the paper is as follows. In Sections \ref{S:monoid} and \ref{S:semiring} we present two approaches to constructing pathwise duality functions. The first approach is based on commutative monoids and the second approach on semirings. In Section~\ref{S:other} we discuss some special cases: a class of duality functions that lie on the intersection of both approaches and duality functions based on lattices that are a special case of the first approach. In Section~\ref{S:numeric} we use computer assisted calculations to find all duality functions that our two approaches yield for local state spaces with cardinality at most four. This includes both well known duality functions and new examples. Although the proofs of our results are quite short, for readability, we have moved them all to Section~\ref{S:proofs}.


\section{Dualities based on commutative monoids}\label{S:monoid}

\subsection{Commutative monoids}

By definition, a \emph{semigroup} is a pair $(S,+)$ where $S$ is a set and $+$ is an associative operation on $S$, i.e.,
\begin{enumerate}
\item $(x+y)+z=x+(y+z)$ $(x,y,z\in S)$.
\end{enumerate}
A semigroup is \emph{commutative} if moreover
\begin{enumerate}\addtocounter{enumi}{1}
\item $x+y=y+x$ $(x,y\in S)$.
\end{enumerate}
A \emph{neutral element} of a semigroup $(S,+)$ is an element $0\in S$ such that
\begin{enumerate}\addtocounter{enumi}{2}
\item $x+0=x=0+x$ $(x\in S)$.
\end{enumerate}
It is easy to see that the neutral element, if it exists, is unique. By definition, a \emph{monoid} is a semigroup $(S,+)$ that is equipped with a neutral element $0$.

If $(S,+)$ and $(T,+)$ are monoids, then a \emph{homomorphism} from $S$ to $T$ is a function $h:S\to T$ such that
\begin{enumerate}
\item $h(x+y)=h(x)+h(y)$ $(x,y\in S)$,
\item $h(0)=0$.
\end{enumerate}
We denote the set of all homomorphisms from $S$ to $T$ by $\Hi(S,T)$. If $h\in\Hi(S,T)$ is a bijection, then it is easy to see that $h^{-1}\in\Hi(T,S)$. In this case, $h$ is called an \emph{isomorphism}. A subset $S'\sub S$ that contains $0$ and is closed under addition is called a \emph{sub-monoid} of $S$. Then $(S',+)$ is itself a monoid with neutral element $0$.

If $(S,+)$ is a semigroup and $\La$ is a set, then we can naturally equip the space $S^\La$ of functions $f:\La\to S$ with the structure of a semigroup by setting
\be
(g+h)(i):=g(i)+h(i)\qquad(g,h\in S^\La,\ i\in\La).
\ee
If $S$ is commutative, then so is $S^\La$, and if $S$ has a neutral element $0$, then $\un{0}$, defined as
\be
\un{0}(i):=0\qquad(i\in\La)
\ee
is the neutral element of $S^\La$. The following simple lemma shows that if $T$ is commutative, then $\Hi(S,T)$ naturally has the structure of a commutative monoid. We call $\Hi(S,T)$ the \emph{$T$-adjoint} of the monoid $S$.

\bl[Adjoint of a monoid]
Let\label{L:adjoint} $S$ and $T$ be monoids and assume that $T$ is commutative. Then $\Hi(S,T)$ is a sub-monoid of $T^S$.
\el

Let $S,T$ be commutative monoids, let $S':=\Hi(S,T)$ denote the $T$-adjoint of $S$ and let $S'':=\Hi(S',T)$ denote the $T$-adjoint of the $T$-adjoint. We claim that there exists a natural homomorphism from $S$ to $S''$. To see this, for each $x\in S$, we define $L_x:\Hi(S,T)\to T$ by
\be\label{Lx}
L_x(h):=h(x)\qquad\big(x\in S,\ h\in\Hi(S,T)\big).
\ee
With this definition, the following lemma holds.

\bl[Adjoint of the adjoint]
Let\label{L:adjadj} $S$ and $T$ be commutative monoids and let $S':=\Hi(S,T)$ and $S'':=\Hi(S',T)$. Then the map $x\mapsto L_x$ is a homomorphism from $S$ to $S''$.
\el

\subsection{Duality of commutative monoids}\label{S:monoidual}

We are now ready for the central definition of this section. Let $R,S$, and $T$ be commutative monoids and let $\psi:S\times R\to T$ be a function. We say that $S$ is \emph{$T$-dual} to $R$ with \emph{duality function} $\psi$ if the following conditions are satisfied:
\begin{enumerate}
\item $\psi(x_1,y)=\psi(x_2,y)$ for all $y\in R$ implies $x_1=x_2$ $(x_1,x_2\in S)$,
\item $\Hi(S,T)=\{\psi(\,\cdot\,,y):y\in R\}$,
\item $\psi(x,y_1)=\psi(x,y_2)$ for all $x\in S$ implies $y_1=y_2$ $(y_1,y_2\in R)$,
\item $\Hi(R,T)=\{\psi(x,\,\cdot\,):x\in S\}$.
\end{enumerate}
Let $S':=\Hi(S,T)$ be the $T$-adjoint of $S$ and let $S'':=\Hi(S',T)$ be the $T$-adjoint of the $T$-adjoint. Borrowing terminology from the theory of Banach spaces, by definition, we say that $S$ is \emph{$T$-reflexive} if the map $x\mapsto L_x$ defined in (\ref{Lx}) is a bijection (and hence an isomorphism) from $S$ to $S''$. The following proposition links duality in the sense we have just defined to the concept of the $T$-adjoint defined in the previous subsection.

\bp[Monoid duality]
Let\label{P:psiprop} $S,R$, and $T$ be commutative monoids and let $S':=\Hi(S,T)$ and $R':=\Hi(R,T)$ be the $T$-adjoints of $S$ and $R$. Then:
\begin{itemize}
\item[{\bf(a)}] If $S$ is $T$-dual to $R$ with duality function $\psi$, then the map $y\mapsto\psi(\,\cdot\,,y)$ is an isomorphism from $R$ to $S'$ and the map $x\mapsto\psi(x,\,\cdot\,)$ is an isomorphism from $S$ to $R'$. Moreover, $S$ and $R$ are $T$-reflexive. 
\item[{\bf(b)}] If $S$ is $T$-reflexive, then $S$ is $T$-dual to $S'$ with duality function
\be\label{psimoi}
\psi(x,h):=h(x)\qquad(x\in S,\ h\in S').
\ee
\end{itemize}
\ep

In Subsection~\ref{S:examp}, we will list all duality functions between monoids of cardinality at most four. Our examples suggest that such duality functions are not rare. The following proposition links the duality functions of Proposition~\ref{P:psiprop} to the concept of a dual map as defined in (\ref{mapdual}).


\bp[Maps having a dual]
Let\label{P:mondu} $S,R$, and $T$ be commutative monoids such that $S$ is $T$-dual to $R$ with duality function $\psi$. Then a map $m:S\to S$ has a dual map $\hat m:R\to R$ with respect to $\psi$ if and only if $m\in\Hi(S,S)$. The dual map $\hat m$, if it exists, is unique and satisfies $\hat m\in\Hi(R,R)$.
\ep

\subsection{Product spaces}

In view of the applications of our results in the theory of interacting particle systems, it is important to pay special attention to product spaces. We have already seen that if $(S,+)$ is a commutative monoid with neutral element $0$ and $\La$ is a set, then the product space $S^\La$ has the structure of a commutative monoid with neutral element $\un 0$. Similarly, if $S_1,\ldots,S_n$ are commutative monoids, then we can naturally equip the product space $S_1\times\cdots\times S_n$ with the structure of a commutative monoid.

We claim that if $S_1,\ldots,S_n$ and $T$ are commutative monoids, then there exists a natural isomorphism $\Hi(S_1,T)\times\cdots\times\Hi(S_n,T)\cong\Hi(S_1\times\cdots\times S_n,T)$. To see this, for each $\fb=(\fb_1,\ldots,\fb_n)\in\Hi(S_1,T)\times\cdots\times\Hi(S_n,T)$, we define a function $F_\fb:S_1\times\cdots\times S_n\to T$ by
\be\label{Ffb}
F_\fb(\xb):=\sum_{i=1}^n\fb_i(\xb_i)\qquad\big(\xb=(\xb_i)_{i\in\La}\in S_1\times\cdots\times S_n\big).
\ee
Here the outcome of the sum does not depend on the summation order, since $T$ is commutative.

\bl[Adjoints of product spaces]
Let\label{L:produa} $S_1,\ldots,S_n$ and $T$ be commutative monoids. Then the map $\fb\mapsto F_\fb$ is an isomorphism from $\Hi(S_1,T)\times\cdots\times\Hi(S_n,T)$ to $\Hi(S_1\times\cdots\times S_n,T)$.
\el

As a simple application of Lemma~\ref{L:produa}, we obtain a characterisation of $\Hi(S^\La,S^\La)$, or somewhat more generally, the set of homomorphisms between two product monoids $S^\La$ and $R^\De$.

\bl[Homomorphisms between product spaces]
Let\label{L:prodhom} $S,R$ be monoids, let $\La,\De$ be finite sets, and let $m:S^\La\to R^\De$ be a map, with $m(x)=\big(m_j(x)\big)_{j\in\De}$. Then one has $m\in\Hi(S^\La,R^\De)$ if and only if there exists a matrix $M=(M_{ij})_{i\in\La,\ j\in\De}$ with $M_{ij}\in\Hi(S,R)$ for each $i\in\La$ and $j\in\De$, such that
\be
m_j(\xb)=\sum_{i\in\La}M_{ij}(\xb_i)\qquad(\xb\in S^\La,\ j\in\De).
\ee
\el

The following proposition says that any duality between commutative monoids can be ``lifted'' to a duality between product spaces. Note that since $T$ is commutative, the sums in (\ref{psi+psi}) and (\ref{psib}) do not depend on the summation order.

\bp[Duality of product spaces]
Let\label{P:prodrefl} $S_1,\ldots,S_n$, $R_1,\ldots,R_n$, and $T$ be commutative monoids and assume that $S_i$ is $T$-dual to $R_i$ with duality function $\psi_i$ $(1\leq i\leq n)$. Then $S_1\times\cdots\times S_n$ is $T$-dual to $R_1\times\cdots\times R_n$ with duality function
\be\label{psi+psi}
\psib(\xb,\yb):=\sum_{i=1}^n\psi_i(\xb_i,\yb_i)
\qquad\big(\xb\in S_1\times\cdots\times S_n,\ \yb\in R_1\times\cdots\times R_n\big).
\ee
In particular, if $\La$ is a finite set and $S$ is $T$-dual to $R$, then $S^\La$ is $T$-dual to $R^\La$ with duality function
\be\label{psib}
\psib(\xb,\yb):=\sum_{i\in\La}\psi(\xb_i,\yb_i)
\qquad\big(\xb\in S^\La,\ \yb\in R^\La\big).
\ee
\ep

Note that by Proposition~\ref{P:mondu}, a map $m:S^\La\to S^\La$ has a dual with respect to the function $\psib$ defined in (\ref{psib}) if and only if $m\in\Hi(S^\La,S^\La)$. By Lemma~\ref{L:prodhom}, maps $m\in\Hi(S^\La,S^\La)$ are uniquely characterised by a matrix with values in $\Hi(S,S)$.

For example, setting $(S,+):=(\{0,1\},\vee)$, one can check that $S$ is $S$-dual to $S$ with duality function $\psi(x,y):=xy$. Defining $\psib$ as in (\ref{psib}) now yields the additive duality function $\psib_{\rm add}$ from (\ref{adccan}). Similarly, setting $S:=\{0,1\}$ but defining $+$ as addition modulo 2 one can again check that $S$ is $S$-dual to $S$ with duality function $\psi(x,y):=xy$. Defining $\psib$ as in (\ref{psib}) now yields the cancellative duality function $\psib_{\rm canc}$ from (\ref{adccan}). Note that in both these examples, when we identify $S$ as a set with $\{0,1\}$ in the way we have just done, then the ``local'' duality function $\psi:S\times S\to S$ is the same, but the ``global'' duality functions $\psib:S^\La\times S^\La\to S$ are still different since the sum on $S$ is defined differently in each example.

\section{Dualities based on semirings}\label{S:semiring}

By definition, a \emph{semiring} is a triple $(S,+,\cdot)$ such that:
\begin{enumerate}
\item $(S,+)$ is a commutative monoid with neutral element $0$,
\item $(S,\cdot)$ is a monoid with neutral element $1$,
\item $x\cdot 0=0=0\cdot x$ for all $x\in S$,
\item $x\cdot(y+z)=x\cdot y+x\cdot z$ and $(x+y)\cdot z=x\cdot z+y\cdot z$ for all $x,y,z\in S$.
\end{enumerate}
Property~(iv) is called \emph{distributivity}. The semiring $(S,+,\cdot)$ is called \emph{commutative} if the monoid $(S,\cdot)$ is.

Let $(S,+,\cdot)$ be a semiring and let $\La$ be a finite set. Then we can equip the monoid $(S^\La,+)$ with additional structure by defining multiplication by scalars from the left and right as
\be
(x\cdot\yb)(i):=x\cdot\yb(i)
\quand
(\yb\cdot x)(i):=\yb(i)\cdot x
\qquad(x\in S,\ \yb\in S^\La,\ i\in\La).
\ee
One can check that with this definition, $S^\La$ becomes an \emph{$S$-module}. In particular, if $S$ is a field, then $S^\La$ is a \emph{linear space over $S$}. Since we will not need the general concepts of $S$-modules and linear spaces, we omit their definitions. Let $\La$ and $\De$ be finite sets and let $\Fi(S^\La,S^\De)$ be the set of all functions $h:S^\La\to S^\De$. Using the conditions
\begin{enumerate}
\item $h(\xb+\yb)=h(\xb)+h(\yb)$ $(\xb,\yb\in S^\La)$, 
\item $h(x\cdot\yb)=x\cdot h(\yb)$ $(x\in S,\ \yb\in S^\La)$,
\item $h(\yb\cdot x)=h(\yb)\cdot x$ $(x\in S,\ \yb\in S^\La)$,
\end{enumerate}
we define sets of functions by
\bc
\dis\Li(S^\La,S^\De)&:=&\dis\big\{h\in\Fi(S^\La,S^\De):h\mbox{ satisfies (i) and (ii)}\big\},\\[5pt]
\dis\Ri(S^\La,S^\De)&:=&\dis\big\{h\in\Fi(S^\La,S^\De):h\mbox{ satisfies (i) and (iii)}\big\}.
\ec
In other words, $\Li(S^\La,S^\De)$ is the set of homomorphisms $h$ from $S^\La$ to $S^\De$, viewed as a left $S$-modules, and likewise $\Ri(S^\La,S^\De)$ is the set of homomorphisms $h$ from $S^\La$ to $S^\De$, viewed as a right $S$-modules. As before, we let $\Hi(S^\La,S^\De)$ denote the set of all homomorphisms $h$ from the monoid $(S^\La,+)$ into $(S^\De,+)$. Note that setting $x=0$ in (ii) yields $h(\un{0})=h(0\cdot\yb)=0\cdot h(\yb)=\un 0$ so $\Li(S^\La,S^\De)\sub\Hi(S^\La,S^\De)$ and similarly $\Ri(S^\La,S^\De)\sub\Hi(S^\La,S^\De)$. If $S$ is commutative, then $\Li(S^\La,S^\De)=\Ri(S^\La,S^\De)$. In particular, if $S$ is a field, then $\Li(S^\La,S^\De)$ is the space of linear functions $h:S^\La\to S^\De$. The following lemma is similar to Lemma~\ref{L:prodhom}.

\bl[Maps between product spaces]
Let\label{L:prodmap} $(S,+,\cdot)$ be a semiring, let $\La,\De$ be finite sets, and let $m:S^\La\to S^\De$ be a map, with $m(x)=\big(m_j(x)\big)_{j\in\De}$. Then one has $m\in\Li(S^\La,S^\De)$ if and only if there exists a matrix $M=(M_{ij})_{i\in\La,\ j\in\De}$ with $M_{ij}\in\Li(S,S)$ for each $i\in\La$ and $j\in\De$, such that
\be
m_j(\xb)=\sum_{i\in\La}M_{ij}(\xb_i)\qquad(\xb\in S^\La,\ j\in\De).
\ee
\el

We define a function $\psib:S^\La\times S^\La\to S$ by
\be\label{linpsi}
\psib(\xb,\yb):=\sum_{i\in\La}\xb(i)\cdot\yb(i)\qquad(\xb,\yb\in S^\La).
\ee
The following lemma says that this function has properties similar to the duality functions of Subsection~\ref{S:monoidual}.

\bl[Duality function for modules over a semiring]
Let\label{L:linpsi} $S$ be a semiring, let $\La$ be a finite set, and let $\psib:S^\La\times S^\La\to S$ be defined as in (\ref{linpsi}). Then:
\begin{enumerate}
\item $\psib(\xb_1,\yb)=\psib(\xb_2,\yb)$ for all $\yb\in S^\La$ implies $\xb_1=\xb_2$ $(\xb_1,\xb_2\in S^\La)$,
\item $\Li(S^\La,S)=\big\{\psib(\,\cdot\,,\yb):\yb\in S^\La\big\}$,
\item $\psib(\xb,\yb_1)=\psib(\xb,\yb_2)$ for all $\xb\in S^\La$ implies $\yb_1=\yb_2$ $(\yb_1,\yb_2\in S^\La)$,
\item $\Ri(S^\La,S)=\big\{\psib(\xb,\,\cdot\,):\xb\in S^\La\big\}$.
\end{enumerate}
\el

The following proposition is similar to Proposition~\ref{P:mondu}.

\bp[Maps having a dual]
Let\label{P:lindu} $S$ be a semiring and let $\La$ be a finite set. Then a map $m:S^\La\to S^\La$ has a dual map $\hat m:S^\La\to S^\La$ with respect to the function $\psib$ defined in (\ref{linpsi}) if and only if $m\in\Li(S^\La,S^\La)$. The dual map $\hat m$, if it exists, is unique and satisfies $\hat m\in\Ri(S^\La,S^\La)$.
\ep

In the special case that $S=\R$, the duality function in (\ref{linpsi}) is the standard inner product on $\R^\La$ and $\hat m$ is the adjoint of the linear map $m$ with respect to this inner product. Linear duality with this duality function has long been used in the study of linear interacting particle systems; see \cite[Chapter~IX]{Lig85} for an overview. It has already been pointed out in \cite[Section~2.6]{Swa13} that linear systems duality can be generalised to linear spaces over arbitrary fields and that in particular, choosing for $S$ the finite field with two elements, one can view cancellative duality as a special case of linear duality. In fact, applying Proposition~\ref{P:lindu} to the semiring $(\{0,1\},\vee,\cdot)$ we see that additive duality also fits into the general class of dualities discussed in the present section.


\section{Some special cases}\label{S:other}

\subsection{Semirings generated by the unit element}\label{S:onegen}


The duality functions in (\ref{psib}) and (\ref{linpsi}) have a similar form. In the present subsection, we will see that under certain conditions, they coincide. Let $(S,+,\cdot)$ be a semiring. Recall that $1\in S$ denotes the neutral element of the product. We say that 1 \emph{generates} $(S,+)$ if each $x\in S$ with $x\neq 0$ is of the form
\[
x=\underbrace{1+\cdots+1}_{\mbox{$n$ times}}
\]
for some integer $n\geq 1$. If $1$ generates $(S,+)$, then it is easy to see that $(S,+,\cdot)$ must be commutative. If $1$ generates $(S,+)$ and $\La,\De$ are finite sets, then we claim that $\Li(S^\La,S^\De)=\Hi(S^\La,S^\De)$, i.e., each $h\in\Hi(S^\La,S^\De)$ satisfies the defining property~(ii) of $\Li(S^\La,S^\De)$. For $x=0$ this is clear since $h(\un 0)=\un 0$. Otherwise, we can write $x=1+\cdots+1$ and observe that
\be\ba{l}
\dis h(x\cdot\yb)=h\big((1+\cdots+1)\cdot\yb\big)
=h\big(\yb+\cdots+\yb\big)\\[5pt]
\dis\quad=h(\yb)+\cdots+h(\yb)=(1+\cdots+1)\cdot h(\yb)=x\cdot h(\yb).
\ec
The following lemma shows that if $1$ generates $(S,+)$, then the semiring-based duality in the sense of Proposition~\ref{P:lindu} is a special case of monoid duality as defined in Subsection~\ref{S:monoidual}.

\bl[Semirings generated by 1]
Assume\label{L:semiring} that $(S,+,\cdot)$ is a commutative semiring and that $1$ generates $(S,+)$. Then $(S,+)$ is $(S,+)$-dual to $(S,+)$ with duality function $\psi(x,y):=x\cdot y$ $(x,y\in S)$.
\el

\subsection{Lattice duality}

A \emph{lattice} is a partially ordered set $(S,\leq)$ with the property that each $x,y\in S$ have a least upper bound $x\vee y$ and a greatest lower bound $x\wedge y$. Following \cite[Subsection~2.4]{SS18}, we say that a lattice $(S^\ast,\leq)$ is \emph{dual} to $(S,\leq)$ if there exists a bijection $S\ni x\mapsto x^\ast\in S^\ast$ such that $x\leq y$ if and only if $x^\ast\geq y^\ast$ $(x,y\in S)$. Clearly, each lattice has a dual, and the dual is unique up to isomorphism. Each finite lattice has unique minimal and maximal elements. If $(S,\leq)$ is a finite lattice with minimal element $0$, then $(S,\vee)$ is a monoid with neutral element $0$. The following lemma says that the monoids $(S,\vee)$ and $(S^\ast,\vee)$ are dual in the sense defined in Subsection~\ref{S:monoidual}.

\bl[Lattice duality]
Let\label{L:lattice} $T$ denote the monoid $(\{0,1\},\vee)$. Let $(S,\leq)$ be a finite lattice and let $(S^\ast,\leq)$ be its dual lattice. Then $(S,\vee)$ is $T$-dual to $(S^\ast,\vee)$ with duality function
\be\label{lattice}
\psi(x,y):=\left\{\ba{ll}
0\quad&\mbox{if }x\leq y^\ast,\\[5pt]
1\quad&\mbox{otherwise.}
\ea\right.\qquad(x\in S,\ y\in S^\ast).
\ee
\el

Pathwise dualities based on dual lattices were studied in \cite{SS18}. In particular, \cite[Lemma~6]{SS18} is just our Proposition~\ref{P:mondu} restricted to the special setting of Lemma~\ref{L:lattice}. Additive duality is a special case of lattice duality, restricted to lattices of the form $\{0,1\}^\La$. As discussed in \cite[Subsection~3.3]{SS18}, the duality of the two-stage contact process discovered by Krone \cite{Kro99} is based on lattices of the form $\{0,1,2\}^\La$.

\section{Examples and discussion}\label{S:numeric}

\subsection{Monoids with up to four elements}\label{S:examp}

Using the approaches in Sections \ref{S:monoid} and \ref{S:semiring}, one can find duality functions of the form
\be
\psib(\xb,\yb)=\sum_{i\in\La}\psi\big(\xb(i),\yb(i)\big)
\qquad(\xb\in S^\La,\ \yb\in R^\La),
\ee
where $\psi:S\times R\to T$ is a ``local'' duality function and the sum is taken in the commutative monoid $(T,+)$. Combining Propositions \ref{P:mondu} and \ref{P:lindu} with Lemmas \ref{L:prodhom} and \ref{L:prodmap}, one can find all maps $m:S^\La\to S^\La$ that have a dual with respect to $\psib$. As explained in Subsection~\ref{S:path}, interacting particle systems based on these maps then have a pathwise dual.

In this section, we will systematically find all local duality functions that arise from these approaches when the spaces $R,S,T$ have cardinality at most four. The number of commutative monoids, up to isomorphism, with $1,2,3,4,5,6,7,\ldots$ elements is $1,2,5,19,78,421,2637,\ldots$ (sequence A058131 in \cite{OEIS}), so beyond cardinality four the sort of brute force approach outlined in Section~\ref{S:moex} quickly becomes impractical.

In the present subsection, we start by listing all commutative monoids with at most four elements. For those with precisely four elements, we have used \cite{For55} as our source. For a monoid of cardinality $n$, we have enumerated its elements $0,\ldots,n-1$ where $0$ always denotes the neutral element. To enumerate the other elements we have applied following rules.
\begin{itemize}
\item[1)] If a monoid $S$ denotes the addition of commutative semirings, and the neutral element of the multiplication is the same one in all those semirings, then we have denoted it by~1.
\item[2)] If a monoid $S$ possesses an \emph{absorbing element} (i.e.\ an element $x\in M$ such that $x+y=y+x=x$ for all $y\in S$), we have denoted it by $n-1$.
\item[3)] If a monoid $S$ possesses an \emph{almost absorbing element} (i.e.\ an element $x\in M$ such that $x+y=y+x=x$ for $y\neq x$ but $x+x\neq x$), we have denoted it by $n-1$.
\end{itemize}
The remaining elements we have denoted in such a way that they appear increasingly often in the addition table. Note that rules 2) and 3) can never contradict themselves. We see, however, several conflicts between rules 1) and 3), where we then have applied rule 1) as indicated by the order.

We have named the monoids $M_0,\ldots,M_{26}$, where $M_0$ is the one monoid with 1 element, $M_1$ and $M_2$ are the monoids with 2 elements, $M_3,\ldots,M_7$ are the ones with 3 elements and $M_8,\ldots,M_{26}$ are the ones with 4 elements. Within the group of monoids with $n$ elements we have ordered the monoids such that the first ones have an absorbing element, then next ones have an almost absorbing element and the ones without either one form the last group. Within these groups we have ordered the monoids such that the $n-1$ appears decreasingly often in the addition table. If multiple monoids within a group have the same number of $(n-1)$-entries in their addition table they are sorted in such a way that the number of $(n-2)$-entries decreases etc.

Below we list the addition tables of $M_0,\ldots,M_7$. The addition tables of $M_8,\ldots,M_{26}$ are given in Appendix~\ref{A:addtable}.

\begin{center}
\begin{tabular}{c|c}
$M_0$&0\\
\hline
0&0\\
\end{tabular}\\[1em]
\begin{tabular}{c|cc}
$M_1$&0&1\\
\hline
0&0&1\\
1&1&1\\
\end{tabular}
\hspace{1em}
\begin{tabular}{c|cc}
$M_2$&0&1\\
\hline
0&0&1\\
1&1&0\\
\end{tabular}\\[1em]
\begin{tabular}{c|ccc}
$M_3$&0&1&2\\
\hline
0&0&1&2\\
1&1&2&2\\
2&2&2&2\\
\end{tabular}
\hspace{1em}
\begin{tabular}{c|ccc}
$M_4$&0&1&2\\
\hline
0&0&1&2\\
1&1&1&2\\
2&2&2&2\\
\end{tabular}
\hspace{1em}
\begin{tabular}{c|ccc}
$M_5$&0&1&2\\
\hline
0&0&1&2\\
1&1&0&2\\
2&2&2&2\\
\end{tabular}\\[0.5em]
\begin{tabular}{c|ccc}
$M_6$&0&1&2\\
\hline
0&0&1&2\\
1&1&2&1\\
2&2&1&2\\
\end{tabular}
\hspace{1em}
\begin{tabular}{c|ccc}
$M_7$&0&1&2\\
\hline
0&0&1&2\\
1&1&2&0\\
2&2&0&1\\
\end{tabular}
\end{center}

\subsection{Semirings with up to four elements}

In Section~\ref{S:semiring}, we studied local duality functions of the form $\psi(x,y)=x\cdot y$ $(x,y\in S)$ where $(S,+,\cdot)$ is a semiring. In the present subsection, we find all local duality functions of this form when $S$ has cardinality between two and four.


Recall that if $(S,+,\cdot)$ is a semiring, then $(S,+)$ is a commutative monoid and $(S,\cdot)$ is a monoid. The monoid $(S,\cdot)$ has an absorbing element, which is the neutral element $0$ of $(S,+)$. It turns out that all monoids with two or three elements that contain an absorbing element are commutative, but there exist two monoids with four elements that contain an absorbing element and are non-commutative. We have named these $N_1$ and $N_2$. Their multiplication tables appear in Appendix~\ref{A:SRmult}. Using a computer, we have found all pairs of monoids $(S,R)$ so that $(S,+)$ is commutative, $(R,\cdot)$ contains an absorbing element, $S$ and $R$ have the same cardinality, which is at most four, and it is possible to identify the elements of $S$ and $R$ in such a way that $(S,+,\cdot)$ is a semiring.

Below we list all possible ways to define a multiplication $\cdot$ on the commutative monoids $M_k$ with $k=1,\ldots,7$ such that the $(M_k,+,\cdot)$ is a semiring. Below each multiplication table, we have indicated to which monoid $(M_k,\cdot)$ is isomorphic. Note that each multiplication table gives rise to a duality function of the form (\ref{linpsi}). The corresponding tables for the monoids $M_8,\ldots,M_{26}$ are given in Appendix~\ref{A:SRmult}. We have only listed semirings that are not isomorphic to each other. In other words, on some of the monoids it may be possible to define a multiplication in a way that is not listed, but in such a case the resulting semiring is isomorphic to a semiring that occurs in our list.

\begin{center}
\begin{tabular}{r|cc}
$(M_1,\cdot)$&0&1\\
\hline
0&0&0\\
1&0&1\\
\multicolumn{3}{l}{\rule{0em}{1.5em}mult. $\cong M_1$}
\end{tabular}
\quad
\begin{tabular}{r|cc}
$(M_2,\cdot)$&0&1\\
\hline
0&0&0\\
1&0&1\\
\multicolumn{3}{l}{\rule{0em}{1.5em}mult. $\cong M_1$}
\end{tabular}\\[2em]

\begin{tabular}{r|ccc}
$(M_3,\cdot)$&0&1&2\\
\hline
0&0&0&0\\
1&0&1&2\\
2&0&2&2\\
\multicolumn{4}{l}{\rule{0em}{1.5em}mult. $\cong M_4$}
\end{tabular}\quad
\begin{tabular}{r|ccc}
$(M_4,\cdot)$&0&1&2\\
\hline
0&0&0&0\\
1&0&0&1\\
2&0&1&2\\
\multicolumn{4}{l}{\rule{0em}{1.5em}mult. $\cong M_3$}
\end{tabular}\quad
\begin{tabular}{r|ccc}
$(M_4,\cdot)$&0&1&2\\
\hline
0&0&0&0\\
1&0&1&2\\
2&0&2&2\\
\multicolumn{4}{l}{\rule{0em}{1.5em}mult. $\cong M_4$}
\end{tabular}\quad
\begin{tabular}{r|ccc}
$(M_4,\cdot)$&0&1&2\\
\hline
0&0&0&0\\
1&0&1&1\\
2&0&1&2\\
\multicolumn{4}{l}{\rule{0em}{1.5em}mult. $\cong M_4$}
\end{tabular}\\[1em]
\begin{tabular}{r|ccc}
$(M_6,\cdot)$&0&1&2\\
\hline
0&0&0&0\\
1&0&1&2\\
2&0&2&2\\
\multicolumn{4}{l}{\rule{0em}{1.5em}mult. $\cong M_4$}
\end{tabular}\quad
\begin{tabular}{r|ccc}
$(M_7,\cdot)$&0&1&2\\
\hline
0&0&0&0\\
1&0&1&2\\
2&0&2&1\\
\multicolumn{4}{l}{\rule{0em}{1.5em}mult. $\cong M_5$}
\end{tabular}
\end{center}

\subsection{Dualities between commutative monoids}\label{S:moex}

We have used a computer to find all quadruples $(R,S,T,\psi)$ such that $R,S,T$ are commutative monoids with cardinality at most four and $S$ is $T$-dual to $R$ with duality function $\psi$, in the sense defined in Subsection~\ref{S:monoidual}. We have proceeded as follows. For each pair $(S,T)$ of commutative monoids with at least two elements each, we used a computer to calculate $\Hi(S,T)$ by brute force, by checking for every function from $S$ to $T$ whether it is a homomorphism. In all cases where $\Hi(S,T)$ has at most four elements, we used a computer to calculate its addition table and find the commutative monoid from our list that it is isomorphic to. The result of this is a table of size $26\times 26$ that lists for each pair $(S,T)$ the monoid $R$ such that $\Hi(S,T)\cong R$, if $R\in\{M_0,\ldots,M_{26}\}$. Using this table, we found all triples $(R,S,T)$ of monoids of cardinality at most four such that $R\cong\Hi(S,T)$ and $S\cong\Hi(R,T)$.

For each such triple $(R,S,T)$ and for each isomorphism $R\ni y\mapsto f_y\in\Hi(S,T)$, we then calculated the function $\psi:S\times R\to T$ defined as
\be\label{candidate}
\psi(x,y):=f_y(x)\qquad(x\in S,\ y\in R).
\ee
By Proposition~\ref{P:psiprop}, $\psi$ is a duality function if and only if $S$ is $T$-reflexive, and each duality function arises in this way. To check that $S$ is $T$-reflexive, we need to check that the map $x\mapsto L_x$ defined in (\ref{Lx}) is a bijection (and hence an isomorphism) from $S$ to $S''$. Equivalently, setting $L'_x(y):=f_y(x)$ $(x\in S,\ y\in R)$, this says that the map $x\mapsto L'_x\in\Hi(R,T)$ is a bijection. In other words, the function in (\ref{candidate}) is a duality function if and only if $S\ni x\mapsto\psi(x,\,\cdot\,)\in\Hi(R,T)$ is a bijection. Since $S\cong\Hi(R,T)$, the sets $S$ and $\Hi(R,T)$ have the same cardinality, so the function in (\ref{candidate}) is a duality function if and only if the functions $\psi(x,\,\cdot\,)$ with $x$ ranging through $S$ are all different from each other.\footnote{Note that the functions $\psi(\,\cdot\,,y)$ with $y$ ranging through $R$ are trivially all different from each other, since $R\ni y\mapsto f_y\in\Hi(S,T)$ is an isomorphism.}

For all triples $(R,S,T)$ of monoids of cardinality at least two and at most four such that $R\cong\Hi(S,T)$ and $S\cong\Hi(R,T)$, and for all choices of the isomorphism $R\ni y\mapsto f_y\in\Hi(S,T)$, we observed that $R$ and $S$ have the same cardinality and that (\ref{candidate}) defines a duality function. In total, in this way, we identified all 110 quadruples $(R,S,T,\psi)$ such that $R,S,T$ are commutative monoids with cardinality at least two and at most four and $S$ is $T$-dual to $R$ with duality function $\psi$.

A lot of these 110 duality functions are trivially related to each other. We will use the following reductions to restrict the number of duality functions and then list only those that are ``essentially'' different.
\begin{itemize}
\item In many of the 110 examples we have found, it turns out that $T$ contains a smaller sub-monoid $\ti T$ so that the duality function $\psi$ takes values in $\ti T$. For this reason, we will only list examples that are \emph{minimal} in the sense that the function values $\{\psi(x,y):x\in S,\ y\in R\}$ generate the monoid $T$.
\item If $S$ is $T$-dual to $R$ with duality function $\psi$ and $R\ni y\mapsto y'\in R$ is an isomorphism, then $S$ is also $T$-dual to $R$ with the duality function $\psi'$ defined as $\psi'(x,y):=\psi(x,y')$ $(x\in S,\ y\in R)$. If several duality functions are related in this way, then we will list only one of them.
\item If $S$ is $T$-dual to $R$ with duality function $\psi$, then $R$ is $T$-dual to $S$ with duality function $\psi^\dgg$ defined as $\psi^\dgg(y,x):=\psi(x,y)$ $(x\in S,\ y\in R)$. If two duality functions are related in this way, then we will list only one of them.
\end{itemize}
After these reductions, we end up with 22 duality functions that are ``essentially'' different. In all examples that are minimal in the sense defined above, we observed that the cardinality of $T$ is not larger than the cardinalities of $R$ and $S$. The following table lists all duality functions $\psi:S\times R\to T$ where $S,R$ have cardinality two or three and $|T|\leq 3$. Those with $|S|=|R|=4$ and $|T|\leq 4$ are listed in Appendix~\ref{A:psi}. Note that the functions listed in these tables are ``local'' duality functions that then give rise to a ``global'' duality function of the form (\ref{psib}).

\begin{center}
\begin{tabular}{r|cc}
\diagbox[height = 2em]{$M_1$}{$M_1$}&0&1\\
\hline
0&0&0\\
1&0&1\\
\multicolumn{3}{l}{\rule{0em}{1.5em}$\psi_1:M_1\times M_1\to M_1$}
\end{tabular}
\quad
\begin{tabular}{r|cc}
\diagbox[height = 2em]{$M_2$}{$M_2$}&0&1\\
\hline
0&0&0\\
1&0&1\\
\multicolumn{3}{l}{\rule{0em}{1.5em}$\psi_2:M_2\times M_2\to M_2$}
\end{tabular}\\[2em]

\begin{tabular}{r|ccc}
\diagbox[height = 2em]{$M_3$}{$M_3$}&0&1&2\\
\hline
0&0&0&0\\
1&0&1&2\\
2&0&2&2\\
\multicolumn{4}{l}{\rule{0em}{1.5em}$\psi_3:M_3\times M_3\to M_3$}
\end{tabular}\quad
\begin{tabular}{r|ccc}
\diagbox[height = 2em]{$M_4$}{$M_4$}&0&1&2\\
\hline
0&0&0&0\\
1&0&0&1\\
2&0&1&1\\
\multicolumn{4}{l}{\rule{0em}{1.5em}$\psi_4:M_4\times M_4\to M_1$}
\end{tabular}\quad
\begin{tabular}{r|ccc}
\diagbox[height = 2em]{$M_5$}{$M_6$}&0&1&2\\
\hline
0&0&0&0\\
1&0&1&0\\
2&0&2&2\\
\multicolumn{4}{l}{\rule{0em}{1.5em}$\psi_5:M_5\times M_6\to M_5$}
\end{tabular}\\[2em]

\begin{tabular}{r|ccc}
\diagbox[height = 2em]{$M_6$}{$M_6$}&0&1&2\\
\hline
0&0&0&0\\
1&0&1&2\\
2&0&2&2\\
\multicolumn{4}{l}{\rule{0em}{1.5em}$\psi_6:M_6\times M_6\to M_6$}
\end{tabular}\quad
\begin{tabular}{r|ccc}
\diagbox[height = 2em]{$M_7$}{$M_7$}&0&1&2\\
\hline
0&0&0&0\\
1&0&1&2\\
2&0&2&1\\
\multicolumn{4}{l}{\rule{0em}{1.5em}$\psi_7:M_7\times M_7\to M_7$}
\end{tabular}
\end{center}

\subsection{Discussion}

We have described two ways to construct pathwise duality functions for interacting particle systems. The first method is based on duality of commutative monoids as described in Subsection~\ref{S:monoidual} and the second method is based on semirings as described in Section~\ref{S:semiring}. As explained in Subsection~\ref{S:onegen}, the two methods partially overlap. By Lemma~\ref{L:semiring}, if $(S,+,\cdot)$ is a semiring in which 1 generates $(S,+)$, then $\psi(x,y):=x\cdot y$ is a duality function in the sense of Subsection~\ref{S:monoidual}. The duality functions $\psi_1,\psi_2,\psi_3,\psi_6,\psi_7,\psi_9,\psi_{22},\psi_{24}$, and $\psi_{26}$ are of this special form and hence occur also in our tables of multiplications in semirings.

Interestingly, we have found one more duality function between commutative monoids that also occurs in our tables of multiplications in semirings. This is $\psi_{23}$, which also occurs in Appendix~\ref{A:SRmult} as the multiplication on $M_{23}$ that is isomorphic to $M_{11}$. In this example, the neutral element of $(M_{23},\cdot)$ does not generate $(M_{23},+)\cong M_1\times M_2$. Nevertheless, one can check that $\Li(M_{23},M_{23})=\Ri(M_{23},M_{23})=\Hi(M_{23},M_{23})$ and hence by Lemmas \ref{L:prodhom} and \ref{L:prodmap} an analogue statement holds for product spaces.

The cyclic groups $C_2,C_3$ and $C_4$ are given by $M_2,M_7$ and $M_{26}$, respectively, i.e.\ always by the last monoid in the group of monoids with $n$ elements. The duality functions $\psi_2,\psi_7$ and $\psi_{26}$ correspond to multiplication modulo $n$. As we already mentioned, they belong to the duality functions of the special form described by Lemma~\ref{L:semiring}. This follows from the fact that $C_n$, equipped with multiplication modulo $n$, is a semiring that is additively generated by 1.

The cyclic groups $C_2$ and $C_3$, equipped with multiplication modulo 2 and 3, respectively, are in fact finite fields. The finite field $(\F_4,+,\cdot)$ with four elements satisfies $(\F_4,+)\cong M_{25}\cong M_2\times M_2$ and $(\F_4,\cdot)\cong M_{18}$. Its multiplication table can be found in Appendix~\ref{A:SRmult}. The unit element of $(\F_4,\cdot)$ does not generate $(\F_4,+)$ and in fact there exist 12 functions from $\F_4$ to itself that are homomorphisms for $(\F_4,+)$ but not elements of $\Li(\F_4,\F_4)=\Ri(\F_4,\F_4)$. Therefore, the multiplication in $\F_4$ does not correspond to a duality between commutative monoids in the sense of Subsection~\ref{S:monoidual}.

The four lattices with 2--4 elements are $M_1\cong(\{0,1\},\vee)$, $M_4\cong(\{0,1,2\},\vee)$, $M_{11}\cong M_1\times M_1$, and $M_{15}\cong(\{0,1,2,3\},\vee)$. Their corresponding duality functions $\psi_1,\psi_4,\psi_{11}$ and $\psi_{15}$ are hence of the form described in Lemma~\ref{L:lattice}. Note that these duality functions are also the only ones that map into $M_1$.

Since $M_k$ is $M_k$-dual to $M_k$ $(k=1,2)$, Proposition~\ref{P:prodrefl} tells us that $M_{11}\cong M_1\times M_1$ is $M_1$-dual to $M_{11}\cong M_1\times M_1$ and that $M_{25}\cong M_2\times M_2$ is $M_2$-dual to $M_{25}\cong M_2\times M_2$. The corresponding duality functions are $\psi_{11}$ and $\psi_{25}$. Since $M_1$ and $M_2$ are naturally submonoids of $M_{23}\cong M_1\times M_2$, the fact that $M_k$ is $M_k$-dual to $M_k$ $(k=1,2)$ trivially implies that $M_k$ is $M_{23}$-dual to $M_k$ $(k=1,2)$ and hence by Proposition~\ref{P:prodrefl} $M_{23}\cong M_1\times M_2$ is $M_{23}$-dual to $M_{23}\cong M_1\times M_2$. It is easy to check that $M_1$ and $M_2$ are also both sub-monoids of $M_5$, so by the same argument $M_{23}$ is also $M_5$-dual to $M_{23}$. The duality functions in these last two cases are $\psi_{23}$ and $\psi_{235}$. We already encountered $\psi_{11}$ and $\psi_{23}$ before since $M_{11}$ is a lattice and since $\psi_{23}$ also occurs in the tables of multiplicative semirings.

If we discard all duality functions between commutative monoids that we have discussed so far, then we are left with the duality function $\psi_5$ from Subsection~\ref{S:moex} and the duality functions $\psi_{10},\psi_{13},\psi_{16},\psi_{17},\psi_{18},\psi_{21}$ from Appendix~\ref{A:psi} that do not have an easy ``explanation''. Of these, $\psi_5,\psi_{13},\psi_{16}$, and $\psi_{18}$ are dualities between different monoids. These duality functions map into $M_5,M_3,M_5$, and $M_{18}$, respectively. The remaining duality functions $\psi_{10},\psi_{17}$ and $\psi_{21}$ are defined on $M_k\times M_k$ with $k=10,17$, and $21$, and map into $M_3,M_5$, and $M_5$ respectively.

Our computer assisted calculations indicate that duality between commutative monoids in the sense of Subsection~\ref{S:monoidual} is not rare, but we are far from a situation where we can classify all examples. It is interesting that in all cases where $R,S,T$ are commutative monoids with cardinality at most four such that $R\cong\mathcal{H}(S,T)$ and $S\cong\mathcal{H}(R,T)$, it turns out that $S$ is $T$-dual to $R$, which is a priori a stronger statement. It is not clear to us if there is a general truth behind this or if there are counterexamples with monoids of larger cardinality.

The approaches for finding duality functions described in Sections \ref{S:monoid} and \ref{S:semiring} have many similarities and in fact partially yield the same duality functions, as can in many cases be understood by applying Lemma~\ref{L:semiring}. It is therefore natural to ask if these two approaches can be unified in an even more general approach. Since we do not see an immediate answer to this question we leave it for further research.

\section{Proofs}\label{S:proofs}

\subsection*{Outline}

In this section, we prove our results. Lemmas \ref{L:adjoint} and \ref{L:adjadj} and Propositions \ref{P:psiprop} and \ref{P:mondu} are proved in Subsection~\ref{S:gener}. Lemmas \ref{L:produa} and \ref{L:prodhom} and Proposition~\ref{P:prodrefl} are proved in Subsection~\ref{S:prodpr}. Lemmas \ref{L:prodmap} and \ref{L:linpsi}, Proposition~\ref{P:lindu}, and Lemma~\ref{L:semiring} are proved in Subsection~\ref{S:smrng}. Lemma~\ref{L:lattice}, finally, is proved in Subsection~\ref{S:latt}.

\subsection{General theory}\label{S:gener}

\bpro[of Lemma~\ref{L:adjoint}]
It is easy to see that $\un{0}\in\Hi(S,T)$, so it remains to show that $f+g\in\Hi(S,T)$ for all $f,g\in\Hi(S,T)$. Indeed, for each $x,y\in S$ and $f,g\in\Hi(S,T)$,
\be\ba{l}
\dis(f+g)(x+y)=f(x+y)+g(x+y)=\big(f(x)+f(y)\big)+\big(g(x)+g(y)\big)\\[5pt]
\dis\quad=\big(f(x)+g(x)\big)+\big(f(y)+g(y)\big)=(f+g)(x)+(f+g)(y),
\ec
where we have used the commutativity of $T$ in the third step. Since moreover $(f+g)(0)=f(0)+g(0)=0+0=0$, this shows that $f+g\in\Hi(S,T)$.
\epro

\bpro[of Lemma~\ref{L:adjadj}]
Since for each $x\in S$,
\be\ba{l}
\dis L_x(f+g)=(f+g)(x)=f(x)+g(x)=L_x(f)+L_x(g),\\[5pt]
\dis L_x(\un{0})=\un{0}(x)=0,
\ec
we see that $L_x$ is a homomorphism from $S'$ to $T$, i.e., $L_x\in S''$. The fact that $x\mapsto L_x$ is a homomorphism from $S$ to $S''$ now follows by writing
\be\ba{l}
\dis L_{x+y}(f)=f(x+y)=f(x)+f(y)=L_x(f)+L_y(f),\\[5pt]
\dis L_0(f)=f(0)=0.
\ec
\epro

\bpro[of Proposition~\ref{P:psiprop}]
Assume that $S$ is $T$-dual to $R$ with duality function $\psi$. Property~(iv) implies that $\psi(x,y_1+y_2)=\psi(x,y_1)+\psi(x,y_2)$ and $\psi(x,0)=0$, so the map $y\mapsto\psi(\,\cdot\,,y)$ is an homomorphism from $R$ to $S'$. By property~(ii), the map $y\mapsto\psi(\,\cdot\,,y)$ is surjective and by property~(i) it is one-to-one, so we conclude that it is an isomorphism. Since $R$ is $T$-dual to $S$ with duality function $\psi^\dgg(y,x):=\psi(x,y)$, the same argument shows that the map $x\mapsto\psi(x,\,\cdot\,)$ is an isomorphism from $S$ to $R'$.

If we identify $R$ with $S'$ using the isomorphism $y\mapsto\psi(\,\cdot\,,y)$, then we can identify the function $L_x:S'\to T$ defined in (\ref{Lx}) with the function $L_x:R\to T$ defined as $L_x(y):=\psi(x,y)$ $(x\in S,\ y\in\R)$. This means that the map $x\mapsto L_x$ from $S$ to $S''$ corresponds to the map $x\mapsto\psi(x,\,\cdot\,)$ from $S$ to $R'$, which we have just shown to be an isomorphism. This proves that $S$ is $T$-reflexive, and by the symmetry between $S$ and $R$, the same is true for $R$.

Assume, conversely, that $S$ is $T$-reflexive. To show that $S$ is $T$-dual to $S'$ with the duality function $\psi$ defined in (\ref{psimoi}), we must show that:
\begin{enumerate}
\item $\psi(x,g)=\psi(x,h)$ for all $x\in S$ implies $g=h$,
\item $\Hi(S,T)=\{\psi(\,\cdot\,,h):h\in S'\}$,
\item $\psi(x,h)=\psi(y,h)$ for all $h\in S'$ implies $x=y$,
\item $\Hi(S',T)=\{\psi(x,\,\cdot\,):x\in S\}$.
\end{enumerate}
Properties (i) and (ii) are trivial consequences of the definition of the adjoint $S'$. By the same argument, if we define $\psi':S'\times S''\to T$ by $\psi'(h,L):=L(h)$ $(h\in S',\ L\in S'')$, then
\begin{enumerate}
\item $\psi(h,L)=\psi'(x,M)$ for all $h\in S'$ implies $L=M$,
\item $\Hi(S',T)=\{\psi'(\,\cdot\,,L):L\in S''\}$.
\end{enumerate}
Since by assumption, $S$ is $T$-reflexive, we may identify $S$ with $S''$. In this identification, we have $\psi'(h,x)=\psi'(h,L_x)=L_x(h)=h(x)=\psi(x,h)$ so properties (i) and (ii) of the function $\psi'$ imply properties (iii) and (iv) of the function $\psi$.
\epro

\bpro[of Proposition~\ref{P:mondu}]
If $m\in\Hi(S,S)$ and $y\in R$, then $x\mapsto\psi\big(m(x),y\big)$ is a homomorphism from $S$ to $T$, so by properties (ii) and (iii) of the definition of a duality function, there exists a unique element $\hat m(y)\in R$ such that $\psi\big(m(x),y\big)=\psi\big(x,\hat m(y)\big)$ for all $x\in S$. This shows that $m$ has a unique dual map $\hat m:R\to R$ with respect to the duality function $\psi$.

Assume, conversely, that $m:S\to S$ has a dual map $\hat m:R\to R$. Then $\psi\big(m(x_1+x_2),y\big)=\psi\big(x_1+x_2,\hat m(y)\big)=\psi\big(x_1,\hat m(y)\big)+\psi\big(x_2,\hat m(y)\big)=\psi\big(m(x_1),y\big)+\psi\big(m(x_2),y\big)=\psi\big(m(x_1)+m(x_2),y\big)$ for all $x_1,x_2\in S$ and $y\in R$, so using property~(i) of a duality function we see that $m(x_1+x_2)=m(x_1)+m(x_2)$ for all $x_1,x_2\in S$. Since moreover $\psi\big(m(0),y\big)=\psi\big(0,\hat m(y)\big)=0$, this proves that $m\in\Hi(S,S)$.

This completes the proof that a map $m:S\to S$ has a dual map $\hat m:R\to R$ with respect to $\psi$ if and only if $m\in\Hi(S,S)$, and moreover shows that such a dual map is unique. Since $\hat m$ has a dual with respect to the duality function $\psi^\dgg(y,x):=\psi(x,y)$, namely, the map $m:S\to S$, by what we have already proved, we must have $\hat m\in\Hi(R,R)$.
\epro

\subsection{Product spaces}\label{S:prodpr}

\bpro[of Lemma~\ref{L:produa}]
We first check that $F_\fb\in\Hi(S_1\times\cdots\times S_n,T)$ for all $\fb\in\Hi(S_1,T)\times\cdots\times\Hi(S_n,T)$. Indeed
\begin{itemize}
\item $\dis F_\fb(\xb+\yb)=\sum_{i=1}^n\fb_i\big((\xb+\yb)_i\big)=\sum_{i=1}^n\fb_i(\xb_i+\yb_i)=\sum_{i=1}^n\big(\fb_i(\xb_i)+\fb_i(\yb_i)\big)\\ =\Big(\sum_{i=1}^n\fb_i(\xb_i)\Big)+\Big(\sum_{i=1}^n\fb_i(\yb_i)\Big)=F_\fb(\xb)+F_\fb(\yb)$,
\item $\dis F_\fb(\un{0})=\sum_{i=1}^n\fb_i(\un{0}_i)=\sum_{i=1}^n\fb_i(0)=\sum_{i=1}^n0=0$.
\end{itemize}
We next check that $\fb\mapsto F_\fb$ is a bijection. We first show that it is one-to-one. For each $1\leq i\leq n$ and $x\in S_i$, let us define $x^i\in S_1\times\cdots\times S_n$ by $x^i_j:=x$ if $i=j$ and $:=0$ otherwise. Then $\fb\neq\gb$ implies $\fb_i\neq\gb_i$ for some $1\leq i\leq n$ and hence there exists an $x\in S_i$ such that $\fb_i(x)\neq\gb_i(x)$. Now $F_\fb(x^i)=\fb_i(x)\neq\gb_i(x)=F_\gb(x^i)$ which shows that $F_\fb\neq F_\gb$. It remains to show that $\fb\mapsto F_\fb$ is surjective. For each $F\in\Hi(S_1\times\cdots\times S_n,T)$, we define $\fb\in\Hi(S_1,T)\times\cdots\times\Hi(S_n,T)$ by $\fb_i(x):=F(x^i)$ $(1\leq i\leq n,\ x\in S_i)$. Then for each $\xb\in S_1\times\cdots\times S_n$, we have
\be
F(\xb)=F\Big(\sum_{i=1}^n(\xb_i)^i\Big)=\sum_{i=1}^nF((\xb_i)^i)=\sum_{i=1}^n\fb_i(\xb_i)=F_\fb(\xb),
\ee
which shows that $F=F_\fb$.

To complete the proof, we must show that $\fb\mapsto F_\fb$ is a homomorphism. We denote the neutral element of $\Hi(S_i,T)$ by $o_i$ and the neutral element of $\Hi(S_1,T)\times\cdots\times\Hi(S_n,T)$ by $\un{o}$. Then 
\begin{itemize}
\item $\dis F_{\fb+\gb}(\xb)=\sum_{i=1}^n(\fb+\gb)_i(\xb_i)=\sum_{i=1}^n(\fb_i+\gb_i)(\xb_i)=\sum_{i=1}^n\big(\fb_i(\xb_i)+\gb_i(\xb_i)\big)=\\ \Big(\sum_{i=1}^n\fb_i(\xb_i)\Big)+\Big(\sum_{i=1}^n\fb_i(\xb_i)\Big)=F_\fb(\xb)+F_\gb(\xb)$,
\item $\dis F_{\un{o}}(\xb)=\sum_{i=1}^n\un{o}_i(\xb_i)=\sum_{i=1}^no_i(\xb_i)=\sum_{i=1}^n0=0$.
\end{itemize}
\epro

\bpro[of Lemma~\ref{L:prodhom}]
This follows from applying Lemma~\ref{L:produa} to the maps $m_j$ for each $j\in\De$.
\epro

\bpro[of Proposition~\ref{P:prodrefl}]
We need to check that $\psib$ satisfies conditions (i)--(iv) of the definition of a duality function. By the symmetry between the $S_i$'s and $R_i$'s, it suffices to check conditions (i) and (ii). Similarly to what we did in the proof of Lemma~\ref{L:produa}, for each $1\leq i\leq n$ and $y\in R_i$, let us define $y^i\in R_1\times\cdots\times R_n$ by $y^i_j:=y$ if $i=j$ and $:=0$ otherwise. Then $\psib(\xb,y^i)=\psi_i(\xb_i,y)$ so $\xb\in S_1\times\cdots\times S_n$ is uniquely determined by the values of $\psib(\xb,y^i)$ for all $1\leq i\leq n$ and $y\in R_i$, proving that $\psib$ satisfies condition~(i). To prove also condition~(ii) we must show that
\be\label{HSLa}
\Hi(S_1\times\cdots\times S_n,T)=\{\psib(\,\cdot\,,\yb):\yb\in R_1\times\cdots\times R_n\}.
\ee
We observe that $\psib(\xb+\xb',\yb)=\psib(\xb,\yb)+\psib(\xb',\yb)$ and $\psib(\un 0,\yb)=0$, which proves the inclusion $\supset$ in (\ref{HSLa}). Conversely, by Lemma~\ref{L:produa}, each $F\in\Hi(S_1\times\cdots\times S_n,T)$ is of the form $F(\xb)=\sum_{i=1}^n\fb_i(\xb_i)$ for some $\fb_i\in\Hi(S_i,T)$ $(1\leq i\leq n)$. Since $S_i$ is $T$-dual to $R_i$ with duality function $\psi_i$, this implies that there exists an $\yb\in R_1\times\cdots\times R_n$ such that $\fb_i=\psi_i(\,\cdot\,,\yb_i)$ for all $1\leq i\leq n$ and hence $F(\xb)=\psib(\xb,\yb)$ for all $\xb\in S_1\times\cdots\times S_n$, proving the inclusion $\sub$ in (\ref{HSLa}).
\epro

\subsection{Semirings}\label{S:smrng}

\bpro[of Lemma~\ref{L:prodmap}]
It suffices to prove the claim when $\De$ consists of a single element. The general statement then follows by applying the more elementary claim to the maps $m_j$ for each $j\in\De$. Thus, we need to show that $m\in\Li(S^\La,S)$ if and only if there exist $(M_i)_{i\in\La}$ with $M_i\in\Li(S,S)$ for each $i\in\La$, such that
\be\label{mlin}
m(\xb)=\sum_{i\in\La}M_i(\xb_i)\qquad(\xb\in S^\La).
\ee
It is straightforward to check that (\ref{mlin}) defines a map $m\in\Li(S^\La,S)$. To see that each element $m\in\Li(S^\La,S)$ is of this form, for each $x\in S$ and $i\in\La$, we define $x^i\in S^\La$ by $x^i_j:=x$ if $i=j$ and $:=0$ otherwise. Given $m\in\Li(S^\La,S)$, we define $M_i:S\to S$ by $M_i(x):=m(x^i)$ $(x\in S,\ i\in\La)$. Then it is straightforward to check that $M_i\in\Li(S,S)$ and $m$ is of the form (\ref{mlin}). Since this is very similar to the proof of Lemma~\ref{L:produa}, we omit the details.
\epro

\bpro[of Lemma~\ref{L:linpsi}]
By symmetry, it suffices to prove properties (i) and (ii). For each $i\in\La$, let $\eb_i\in S^\La$ be defined as $\eb_i(i):=1$ and $\eb_i(j):=0$ for all $j\in\La\beh\{i\}$. Then $\psi(\xb_1,\yb)=\psi(\xb_2,\yb)$ for all $\yb\in S^\La$ implies $\xb_1(i)=\psi(\xb_1,\eb_i)=\psi(\xb_2,\eb_i)=\xb_2(i)$ for all $i\in\La$ and hence $\xb_1=\xb_2$, proving (i). Using the distributive property of the product and the commutativity of the sum, we see that
\be\ba{l}
\dis\psi(\xb_1+\xb_2,\yb)=\sum_{i\in\La}\big(\xb_1(i)+\xb_2(i)\big)\cdot\yb(i)
=\sum_{i\in\La}\big(\xb_1(i)\cdot\yb(i)+\xb_2(i)\cdot\yb(i)\big)\\[5pt]
\dis\quad=\sum_{i\in\La}\xb_1(i)\cdot\yb(i)+\sum_{i\in\La}\xb_2(i)\cdot\yb(i)
=\psi(\xb_1,y)+\psi(\xb_2,\yb)\qquad(\xb_1,\xb_2,\yb\in S^\La).
\ec
Using the associative and distributive properties of the product, we obtain moreover that
\be
\psi(z\cdot\xb,\yb)=\sum_{i\in\La}\big(z\cdot\xb(i)\big)\cdot\yb(i)
=\sum_{i\in\La}z\cdot\big(\xb(i)\cdot\yb(i)\big)
=z\cdot\sum_{i\in\La}\xb(i)\cdot\yb(i)=z\cdot\psi(\xb,\yb)
\ee
$(\xb,\yb\in S^\La,\ z\in S)$. Applying this with $z=0$, using the fact that $0\cdot x=0$ $(x\in S)$, we see that moreover $\psi(\un 0,\yb)=0$ $(\yb\in S^\La)$, so we conclude that $\Li(S^\La,S)\supset\big\{\psi(\,\cdot\,,\yb):\yb\in S^\La\big\}$.

To prove the reverse inclusion, assume that $h\in\Li(S^\La,S)$. We will prove that $h=\psi(\,\cdot\,,\yb)$ with $\yb(i):=h(\eb_i)$ $(i\in\La)$. Indeed,
\be
h(\xb)=h\big(\sum_{i\in\La}\xb(i)\cdot\eb_i\big)=\sum_{i\in\La}\xb(i)\cdot h\big(\eb_i\big)=\psi(\xb,\yb)\qquad(\xb\in S^\La),
\ee
which concludes our proof.
\epro

\bpro[of Proposition~\ref{P:lindu}]
If $m\in\Li(S^\La,S^\La)$ and $\yb\in S^\La$, then by Lemma~\ref{L:linpsi}~(ii), the map $\xb\mapsto\psib\big(m(\xb),\yb\big)$ is an element of $\Li(S^\La,S)$, so by Lemma~\ref{L:linpsi} (i) and (ii), there exists a unique element $\hat m(\yb)\in S^\La$ such that $\psib\big(m(\xb),\yb\big)=\psib\big(\xb,\hat m(\yb)\big)$ for all $\xb\in S^\La$. This shows that $m$ has a unique dual map $\hat m$ with respect to the duality function $\psib$.

Assume, conversely, that $m:S^\La\to S^\La$ has a dual map $\hat m:S^\La\to S^\La$ with respect to the duality function $\psib$. Then by Lemma~\ref{L:linpsi}~(ii),
\begin{itemize}
\item $\psib\big(m(\xb_1+\xb_2),\yb\big)=\psib\big(\xb_1+\xb_2,\hat m(\yb)\big)=\psib\big(\xb_1,\hat m(\yb)\big)+\psib\big(\xb_2,\hat m(\yb)\big)\\ =\psib\big(m(\xb_1),\yb\big)+\psib\big(m(\xb_2),\yb\big)=\psib\big(m(\xb_1)+m(\xb_2),\yb\big)$ $(\xb_1,\xb_2,\yb\in S^\La)$
\item  $\psib\big(m(z\cdot\xb),\yb\big)=\psib\big(z\cdot\xb,\hat m(\yb)\big)=z\cdot\psib\big(\xb,\hat m(\yb)\big)=z\cdot\psib\big(m(\xb),\yb\big)\\ =\psib\big(z\cdot m(\xb),\yb\big)$ $(\xb,\yb\in S^\La,\ z\in S)$.
\end{itemize}
Since this holds for all $\yb\in S^\La$, by Lemma~\ref{L:linpsi}~(i), we conclude that $m\in\Li(S^\La,S^\La)$.

This completes the proof that a map $m:S^\La\to S^\La$ has a dual map $\hat m:S^\La\to S^\La$ with respect to $\psib$ if and only if $m\in\Li(S^\La,S^\La)$, and moreover shows that such a dual map is unique. In exactly the same way, using Lemma~\ref{L:linpsi} (iii) and (iv), we see that a map $\hat n:S^\La\to S^\La$ has a dual map $n:S^\La\to S^\La$ with respect to the duality function $\psib^\dgg(\yb,\xb):=\psib(\xb,\yb)$ $(\xb,\yb\in S^\La)$ if and only if $\hat n\in\Ri(S^\La,S^\La)$. Applying this to $\hat n=\hat m$, which has $m$ as a dual map, we see that $\hat m\in\Ri(S^\La,S^\La)$.
\epro

\bpro[of Lemma~\ref{L:semiring}]
Immediate from Lemma~\ref{L:linpsi} and the observation that $\Li(S^\La,S)=\Ri(S^\La,S)=\Hi(S^\La,S)$.
\epro

\subsection{Lattices}\label{S:latt}

\bpro[of Lemma~\ref{L:lattice}]
Let $\{x\leq y^\ast\}=\{y\leq x^\ast\}$ denote the set of all $(x,y)\in S\times S^\ast$ such that $x\leq y^\ast$, and let $1_{\{x\leq y^\ast\}}$ denote its indicator function. Set $R:=(\{0,1\},\wedge)$ and $\ti\psi(x,y):=1_{\{x\leq y^\ast\}}$ $(x\in S,\ y\in S^\ast)$. Then we may equivalently prove that $S$ is $R$-dual to $S^\ast$ with duality function $\ti\psi$. We check conditions (i)--(iv) of our definition of duality of commutative monoids in Subsection~\ref{S:monoidual}. By symmetry, it suffices to check conditions (i) and (ii). Condition~(i) follows from the fact that $y\mapsto y^\ast$ is a bijection and $1_{\{x_1\leq z\}}=1_{\{x_2\leq z\}}$ for all $z\in S$, and in particular for $z=x_1,x_2$, implies $x_1\leq x_2\leq x_1$ and hence $x_1=x_2$. To check condition~(ii), we first observe that
\be
1_{\{0\leq y^\ast\}}=1\quand
1_{\{x_1\vee x_2\leq y^\ast\}}=1_{\{x_1\leq y^\ast\}}\wedge 1_{\{x_0\leq y^\ast\}}
\qquad(x_1,x_2\in S,\ y\in S^\ast).
\ee
Since $1$ is the neutral element of $R$, this shows that $\psi(\,\cdot\,,y)\in\Hi(S,R)$ for all $y\in S^\ast$. Assume, conversely, that $h\in\Hi(S,R)$. To complete the proof, we must show that $h(x)=1_{\{x\leq z\}}$ $(x\in S)$ for some $z\in S$. Since $h(0)=1$, the set $\{x:h(x)=1\}$ is nonempty, so using the finiteness of $S$ we can define $z:=\bigvee\{x:h(x)=1\}$. We observe that $h(x_1)=1=h(x_2)$ implies
\be
h(x_1\vee x_2)=h(x_1)\wedge h(x_2)=1\wedge 1=1.
\ee
It follows that $h(z)=1$ and more generally $h(x)=h(x\vee z)=h(z)=1$ for all $x\leq z$. Conversely, $h(x)=1$ implies that $x$ is an element of $\{x:h(x)=1\}$ and hence $x\leq z$ by the definition of $z$.
\epro

\appendix

\section{Appendix}

\subsection{Addition tables of commutative monoids of order four}\label{A:addtable}

\begin{center}
\begin{tabular}{r|cccc}
$M_8$&0&1&2&3\\
\hline
0&0&1&2&3\\
1&1&3&3&3\\
2&2&3&3&3\\
3&3&3&3&3\\
\end{tabular}
\hspace{1em}
\begin{tabular}{r|cccc}
$M_9$&0&1&2&3\\
\hline
0&0&1&2&3\\
1&1&2&3&3\\
2&2&3&3&3\\
3&3&3&3&3\\
\end{tabular}
\hspace{1em}
\begin{tabular}{r|cccc}
$M_{10}$&0&1&2&3\\
\hline
0&0&1&2&3\\
1&1&3&3&3\\
2&2&3&2&3\\
3&3&3&3&3\\
\end{tabular}\\[0.5em]
\begin{tabular}{r|cccc}
$M_{11}$&0&1&2&3\\
\hline
0&0&1&2&3\\
1&1&1&3&3\\
2&2&3&2&3\\
3&3&3&3&3\\
\end{tabular}
\hspace{1em}
\begin{tabular}{r|cccc}
$M_{12}$&0&1&2&3\\
\hline
0&0&1&2&3\\
1&1&0&2&3\\
2&2&2&3&3\\
3&3&3&3&3\\
\end{tabular}
\hspace{1em}
\begin{tabular}{r|cccc}
$M_{13}$&0&1&2&3\\
\hline
0&0&1&2&3\\
1&1&3&1&3\\
2&2&1&2&3\\
3&3&3&3&3\\
\end{tabular}\\[0.5em]
\begin{tabular}{r|cccc}
$M_{14}$&0&1&2&3\\
\hline
0&0&1&2&3\\
1&1&2&2&3\\
2&2&2&2&3\\
3&3&3&3&3\\
\end{tabular}
\hspace{1em}
\begin{tabular}{r|cccc}
$M_{15}$&0&1&2&3\\
\hline
0&0&1&2&3\\
1&1&1&2&3\\
2&2&2&2&3\\
3&3&3&3&3\\
\end{tabular}
\hspace{1em}
\begin{tabular}{r|cccc}
$M_{16}$&0&1&2&3\\
\hline
0&0&1&2&3\\
1&1&0&2&3\\
2&2&2&2&3\\
3&3&3&3&3\\
\end{tabular}\\[0.5em]
\begin{tabular}{r|cccc}
$M_{17}$&0&1&2&3\\
\hline
0&0&1&2&3\\
1&1&2&1&3\\
2&2&1&2&3\\
3&3&3&3&3\\
\end{tabular}
\hspace{1em}
\begin{tabular}{r|cccc}
$M_{18}$&0&1&2&3\\
\hline
0&0&1&2&3\\
1&1&2&0&3\\
2&2&0&1&3\\
3&3&3&3&3\\
\end{tabular}
\hspace{1em}
\begin{tabular}{r|cccc}
$M_{19}$&0&1&2&3\\
\hline
0&0&1&2&3\\
1&1&2&2&3\\
2&2&2&2&3\\
3&3&3&3&2\\
\end{tabular}\\[0.5em]
\begin{tabular}{r|cccc}
$M_{20}$&0&1&2&3\\
\hline
0&0&1&2&3\\
1&1&3&1&1\\
2&2&1&2&3\\
3&3&1&3&3\\
\end{tabular}
\hspace{1em}
\begin{tabular}{r|cccc}
$M_{21}$&0&1&2&3\\
\hline
0&0&1&2&3\\
1&1&3&1&1\\
2&2&1&0&3\\
3&3&1&3&3\\
\end{tabular}
\hspace{1em}
\begin{tabular}{r|cccc}
$M_{22}$&0&1&2&3\\
\hline
0&0&1&2&3\\
1&1&3&3&2\\
2&2&3&3&2\\
3&3&2&2&3\\
\end{tabular}\\[0.5em]
\begin{tabular}{r|cccc}
$M_{23}$&0&1&2&3\\
\hline
0&0&1&2&3\\
1&1&3&3&1\\
2&2&3&0&1\\
3&3&1&1&3\\
\end{tabular}
\hspace{1em}
\begin{tabular}{r|cccc}
$M_{24}$&0&1&2&3\\
\hline
0&0&1&2&3\\
1&1&2&3&1\\
2&2&3&1&2\\
3&3&1&2&3\\
\end{tabular}
\hspace{1em}
\begin{tabular}{r|cccc}
$M_{25}$&0&1&2&3\\
\hline
0&0&1&2&3\\
1&1&0&3&2\\
2&2&3&0&1\\
3&3&2&1&0\\
\end{tabular}\\[0.5em]
\begin{tabular}{r|cccc}
$M_{26}$&0&1&2&3\\
\hline
0&0&1&2&3\\
1&1&2&3&0\\
2&2&3&0&1\\
3&3&0&1&2\\
\end{tabular}
\end{center}

\subsection{Duality functions for commutative monoids of order four}\label{A:psi}

\begin{center}
\begin{tabular}{r|cccc}
\diagbox[height = 2em]{$M_9$}{$M_9$}&0&1&2&3\\
\hline
0&0&0&0&0\\
1&0&1&2&3\\
2&0&2&3&3\\
3&0&3&3&3\\
\multicolumn{5}{l}{\rule{0em}{1.5em}$\psi_9:M_9\times M_9\to M_9$}
\end{tabular}\quad
\begin{tabular}{r|cccc}
\diagbox[height = 2em]{$M_{10}$}{$M_{10}$}&0&1&2&3\\
\hline
0&0&0&0&0\\
1&0&1&2&2\\
2&0&2&0&2\\
3&0&2&2&2\\
\multicolumn{5}{l}{\rule{0em}{1.5em}$\psi_{10}:M_{10}\times M_{10}\to M_3$}
\end{tabular}\quad
\begin{tabular}{r|cccc}
\diagbox[height = 2em]{$M_{11}$}{$M_{11}$}&0&1&2&3\\
\hline
0&0&0&0&0\\
1&0&0&1&1\\
2&0&1&0&1\\
3&0&1&1&1\\
\multicolumn{5}{l}{\rule{0em}{1.5em}$\psi_{11}:M_{11}\times M_{11}\to M_1$}
\end{tabular}\\[2em]

\begin{tabular}{r|cccc}
\diagbox[height = 2em]{$M_{13}$}{$M_{14}$}&0&1&2&3\\
\hline
0&0&0&0&0\\
1&0&1&2&2\\
2&0&0&0&2\\
3&0&2&2&2\\
\multicolumn{5}{l}{\rule{0em}{1.5em}$\psi_{13}:M_{13}\times M_{14}\to M_3$}
\end{tabular}\quad
\begin{tabular}{r|cccc}
\diagbox[height = 2em]{$M_{15}$}{$M_{15}$}&0&1&2&3\\
\hline
0&0&0&0&0\\
1&0&0&0&1\\
2&0&0&1&1\\
3&0&1&1&1\\
\multicolumn{5}{l}{\rule{0em}{1.5em}$\psi_{15}:M_{15}\times M_{15}\to M_1$}
\end{tabular}\quad
\begin{tabular}{r|cccc}
\diagbox[height = 2em]{$M_{16}$}{$M_{20}$}&0&1&2&3\\
\hline
0&0&0&0&0\\
1&0&1&0&0\\
2&0&2&0&2\\
3&0&2&2&2\\
\multicolumn{5}{l}{\rule{0em}{1.5em}$\psi_{16}:M_{16}\times M_{20}\to M_5$}
\end{tabular}\\[2em]

\begin{tabular}{r|cccc}
\diagbox[height = 2em]{$M_{17}$}{$M_{17}$}&0&1&2&3\\
\hline
0&0&0&0&0\\
1&0&1&0&2\\
2&0&0&0&2\\
3&0&2&2&2\\
\multicolumn{5}{l}{\rule{0em}{1.5em}$\psi_{17}:M_{17}\times M_{17}\to M_5$}
\end{tabular}\quad
\begin{tabular}{r|cccc}
\diagbox[height = 2em]{$M_{18}$}{$M_{24}$}&0&1&2&3\\
\hline
0&0&0&0&0\\
1&0&1&2&0\\
2&0&2&1&0\\
3&0&3&3&3\\
\multicolumn{5}{l}{\rule{0em}{1.5em}$\psi_{18}:M_{18}\times M_{24}\to M_{18}$}
\end{tabular}\quad
\begin{tabular}{r|cccc}
\diagbox[height = 2em]{$M_{21}$}{$M_{21}$}&0&1&2&3\\
\hline
0&0&0&0&0\\
1&0&2&1&2\\
2&0&1&0&0\\
3&0&2&0&2\\
\multicolumn{5}{l}{\rule{0em}{1.5em}$\psi_{21}:M_{21}\times M_{21}\to M_5$}
\end{tabular}\\[2em]

\begin{tabular}{r|cccc}
\diagbox[height = 2em]{$M_{22}$}{$M_{22}$}&0&1&2&3\\
\hline
0&0&0&0&0\\
1&0&1&2&3\\
2&0&2&2&3\\
3&0&3&3&3\\
\multicolumn{5}{l}{\rule{0em}{1.5em}$\psi_{22}:M_{22}\times M_{22}\to M_{22}$}
\end{tabular}\quad
\begin{tabular}{r|cccc}
\diagbox[height = 2em]{$M_{23}$}{$M_{23}$}&0&1&2&3\\
\hline
0&0&0&0&0\\
1&0&1&2&3\\
2&0&2&2&0\\
3&0&3&0&3\\
\multicolumn{5}{l}{\rule{0em}{1.5em}$\psi_{23}:M_{23}\times M_{23}\to M_{23}$}
\end{tabular}\quad
\begin{tabular}{r|cccc}
\diagbox[height = 2em]{$M_{23}$}{$M_{23}$}&0&1&2&3\\
\hline
0&0&0&0&0\\
1&0&2&1&2\\
2&0&1&1&0\\
3&0&2&0&2\\
\multicolumn{5}{l}{\rule{0em}{1.5em}$\psi_{235}:M_{23}\times M_{23}\to M_5$}
\end{tabular}\\[2em]

\begin{tabular}{r|cccc}
\diagbox[height = 2em]{$M_{24}$}{$M_{24}$}&0&1&2&3\\
\hline
0&0&0&0&0\\
1&0&1&2&3\\
2&0&2&1&3\\
3&0&3&3&3\\
\multicolumn{5}{l}{\rule{0em}{1.5em}$\psi_{24}:M_{24}\times M_{24}\to M_{24}$}
\end{tabular}\quad
\begin{tabular}{r|cccc}
\diagbox[height = 2em]{$M_{25}$}{$M_{25}$}&0&1&2&3\\
\hline
0&0&0&0&0\\
1&0&0&1&1\\
2&0&1&0&1\\
3&0&1&1&0\\
\multicolumn{5}{l}{\rule{0em}{1.5em}$\psi_{25}:M_{25}\times M_{25}\to M_{2}$}
\end{tabular}\quad
\begin{tabular}{r|cccc}
\diagbox[height = 2em]{$M_{26}$}{$M_{26}$}&0&1&2&3\\
\hline
0&0&0&0&0\\
1&0&1&2&3\\
2&0&2&0&2\\
3&0&3&2&1\\
\multicolumn{5}{l}{\rule{0em}{1.5em}$\psi_{26}:M_{26}\times M_{26}\to M_{26}$}
\end{tabular}
\end{center}

\subsection{Multiplications in semirings of cardinality four}\label{A:SRmult}

\begin{center}
\begin{tabular}{r|cccc}
$(M_8,\cdot)$&0&1&2&3\\
\hline
0&0&0&0&0\\
1&0&1&2&3\\
2&0&2&3&3\\
3&0&3&3&3\\
\multicolumn{5}{l}{\rule{0em}{1.5em}mult. $\cong M_{14}$}
\end{tabular}\quad
\begin{tabular}{r|cccc}
$(M_8,\cdot)$&0&1&2&3\\
\hline
0&0&0&0&0\\
1&0&1&2&3\\
2&0&2&2&3\\
3&0&3&3&3\\
\multicolumn{5}{l}{\rule{0em}{1.5em}mult. $\cong M_{15}$}
\end{tabular}\quad
\begin{tabular}{r|cccc}
$(M_8,\cdot)$&0&1&2&3\\
\hline
0&0&0&0&0\\
1&0&1&2&3\\
2&0&2&1&3\\
3&0&3&3&3\\
\multicolumn{5}{l}{\rule{0em}{1.5em}mult. $\cong M_{16}$}
\end{tabular}\\[1em]

\begin{tabular}{r|cccc}
$(M_9,\cdot)$&0&1&2&3\\
\hline
0&0&0&0&0\\
1&0&1&2&3\\
2&0&2&3&3\\
3&0&3&3&3\\
\multicolumn{5}{l}{\rule{0em}{1.5em}mult. $\cong M_{14}$}
\end{tabular}\quad
\begin{tabular}{r|cccc}
$(M_{10},\cdot)$&0&1&2&3\\
\hline
0&0&0&0&0\\
1&0&1&2&3\\
2&0&2&0&2\\
3&0&3&2&3\\
\multicolumn{5}{l}{\rule{0em}{1.5em}mult. $\cong M_{13}$}
\end{tabular}\quad
\begin{tabular}{r|cccc}
$(M_{10},\cdot)$&0&1&2&3\\
\hline
0&0&0&0&0\\
1&0&1&2&3\\
2&0&2&2&2\\
3&0&3&2&3\\
\multicolumn{5}{l}{\rule{0em}{1.5em}mult. $\cong M_{15}$}
\end{tabular}\\[1em]

\begin{tabular}{r|cccc}
$(M_{11},\cdot)$&0&1&2&3\\
\hline
0&0&0&0&0\\
1&0&1&0&1\\
2&0&0&2&2\\
3&0&1&2&3\\
\multicolumn{5}{l}{\rule{0em}{1.5em}mult. $\cong M_{11}$}
\end{tabular}\quad
\begin{tabular}{r|cccc}
$(M_{11},\cdot)$&0&1&2&3\\
\hline
0&0&0&0&0\\
1&0&1&2&3\\
2&0&2&0&2\\
3&0&3&2&3\\
\multicolumn{5}{l}{\rule{0em}{1.5em}mult. $\cong M_{13}$}
\end{tabular}\quad
\begin{tabular}{r|cccc}
$(M_{11},\cdot)$&0&1&2&3\\
\hline
0&0&0&0&0\\
1&0&1&2&3\\
2&0&2&3&3\\
3&0&3&3&3\\
\multicolumn{5}{l}{\rule{0em}{1.5em}mult. $\cong M_{14}$}
\end{tabular}\\[1em]

\begin{tabular}{r|cccc}
$(M_{11},\cdot)$&0&1&2&3\\
\hline
0&0&0&0&0\\
1&0&1&2&3\\
2&0&2&2&2\\
3&0&3&2&3\\
\multicolumn{5}{l}{\rule{0em}{1.5em}mult. $\cong M_{15}$}
\end{tabular}\quad
\begin{tabular}{r|cccc}
$(M_{11},\cdot)$&0&1&2&3\\
\hline
0&0&0&0&0\\
1&0&1&2&3\\
2&0&2&1&3\\
3&0&3&3&3\\
\multicolumn{5}{l}{\rule{0em}{1.5em}mult. $\cong M_{16}$}
\end{tabular}\quad
\begin{tabular}{r|cccc}
$(M_{13},\cdot)$&0&1&2&3\\
\hline
0&0&0&0&0\\
1&0&1&2&3\\
2&0&2&0&2\\
3&0&3&2&3\\
\multicolumn{5}{l}{\rule{0em}{1.5em}mult. $\cong M_{13}$}
\end{tabular}\\[1em]

\begin{tabular}{r|cccc}
$(M_{13},\cdot)$&0&1&2&3\\
\hline
0&0&0&0&0\\
1&0&1&2&3\\
2&0&2&2&2\\
3&0&3&2&3\\
\multicolumn{5}{l}{\rule{0em}{1.5em}mult. $\cong M_{15}$}
\end{tabular}\quad
\begin{tabular}{r|cccc}
$(M_{14},\cdot)$&0&1&2&3\\
\hline
0&0&0&0&0\\
1&0&1&2&3\\
2&0&2&2&3\\
3&0&3&3&3\\
\multicolumn{5}{l}{\rule{0em}{1.5em}mult. $\cong M_{15}$}
\end{tabular}\quad
\begin{tabular}{r|cccc}
$(M_{15},\cdot)$&0&1&2&3\\
\hline
0&0&0&0&0\\
1&0&0&0&1\\
2&0&0&0&2\\
3&0&1&2&3\\
\multicolumn{5}{l}{\rule{0em}{1.5em}mult. $\cong M_{8}$}
\end{tabular}\\[1em]

\begin{tabular}{r|cccc}
$(M_{15},\cdot)$&0&1&2&3\\
\hline
0&0&0&0&0\\
1&0&0&0&1\\
2&0&0&1&2\\
3&0&1&2&3\\
\multicolumn{5}{l}{\rule{0em}{1.5em}mult. $\cong M_{9}$}
\end{tabular}\quad
\begin{tabular}{r|cccc}
$(M_{15},\cdot)$&0&1&2&3\\
\hline
0&0&0&0&0\\
1&0&0&0&1\\
2&0&0&2&2\\
3&0&1&2&3\\
\multicolumn{5}{l}{\rule{0em}{1.5em}mult. $\cong M_{10}$}
\end{tabular}\quad
\begin{tabular}{r|cccc}
$(M_{15},\cdot)$&0&1&2&3\\
\hline
0&0&0&0&0\\
1&0&0&1&1\\
2&0&1&2&3\\
3&0&1&3&3\\
\multicolumn{5}{l}{\rule{0em}{1.5em}mult. $\cong M_{13}$}
\end{tabular}\\[1em]

\begin{tabular}{r|cccc}
$(M_{15},\cdot)$&0&1&2&3\\
\hline
0&0&0&0&0\\
1&0&0&1&1\\
2&0&1&2&2\\
3&0&1&2&3\\
\multicolumn{5}{l}{\rule{0em}{1.5em}mult. $\cong M_{13}$}
\end{tabular}\quad
\begin{tabular}{r|cccc}
$(M_{15},\cdot)$&0&1&2&3\\
\hline
0&0&0&0&0\\
1&0&1&2&3\\
2&0&2&3&3\\
3&0&3&3&3\\
\multicolumn{5}{l}{\rule{0em}{1.5em}mult. $\cong M_{14}$}
\end{tabular}\quad
\begin{tabular}{r|cccc}
$(M_{15},\cdot)$&0&1&2&3\\
\hline
0&0&0&0&0\\
1&0&1&1&1\\
2&0&1&1&2\\
3&0&1&2&3\\
\multicolumn{5}{l}{\rule{0em}{1.5em}mult. $\cong M_{14}$}
\end{tabular}\\[1em]

\begin{tabular}{r|cccc}
$(M_{15},\cdot)$&0&1&2&3\\
\hline
0&0&0&0&0\\
1&0&1&2&3\\
2&0&2&2&3\\
3&0&3&3&3\\
\multicolumn{5}{l}{\rule{0em}{1.5em}mult. $\cong M_{15}$}
\end{tabular}\quad
\begin{tabular}{r|cccc}
$(M_{15},\cdot)$&0&1&2&3\\
\hline
0&0&0&0&0\\
1&0&1&1&3\\
2&0&1&2&3\\
3&0&3&3&3\\
\multicolumn{5}{l}{\rule{0em}{1.5em}mult. $\cong M_{15}$}
\end{tabular}\quad
\begin{tabular}{r|cccc}
$(M_{15},\cdot)$&0&1&2&3\\
\hline
0&0&0&0&0\\
1&0&1&1&1\\
2&0&1&2&3\\
3&0&1&3&3\\
\multicolumn{5}{l}{\rule{0em}{1.5em}mult. $\cong M_{15}$}
\end{tabular}\\[1em]

\begin{tabular}{r|cccc}
$(M_{15},\cdot)$&0&1&2&3\\
\hline
0&0&0&0&0\\
1&0&1&1&1\\
2&0&1&2&2\\
3&0&1&2&3\\
\multicolumn{5}{l}{\rule{0em}{1.5em}mult. $\cong M_{15}$}
\end{tabular}\quad
\begin{tabular}{r|cccc}
$(M_{15},\cdot)$&0&1&2&3\\
\hline
0&0&0&0&0\\
1&0&0&0&1\\
2&0&1&2&2\\
3&0&1&2&3\\
\multicolumn{5}{l}{\rule{0em}{1.5em}mult. $\cong N_{1}$}
\end{tabular}\quad
\begin{tabular}{r|cccc}
$(M_{15},\cdot)$&0&1&2&3\\
\hline
0&0&0&0&0\\
1&0&1&1&1\\
2&0&1&2&3\\
3&0&3&3&3\\
\multicolumn{5}{l}{\rule{0em}{1.5em}mult. $\cong N_{2}$}
\end{tabular}\\[1em]

\begin{tabular}{r|cccc}
$(M_{17},\cdot)$&0&1&2&3\\
\hline
0&0&0&0&0\\
1&0&1&2&3\\
2&0&2&2&3\\
3&0&3&3&3\\
\multicolumn{5}{l}{\rule{0em}{1.5em}mult. $\cong M_{15}$}
\end{tabular}\quad
\begin{tabular}{r|cccc}
$(M_{20},\cdot)$&0&1&2&3\\
\hline
0&0&0&0&0\\
1&0&1&2&3\\
2&0&2&0&2\\
3&0&3&2&3\\
\multicolumn{5}{l}{\rule{0em}{1.5em}mult. $\cong M_{13}$}
\end{tabular}\quad
\begin{tabular}{r|cccc}
$(M_{20},\cdot)$&0&1&2&3\\
\hline
0&0&0&0&0\\
1&0&1&2&3\\
2&0&2&2&2\\
3&0&3&2&3\\
\multicolumn{5}{l}{\rule{0em}{1.5em}mult. $\cong M_{15}$}
\end{tabular}\\[1em]

\begin{tabular}{r|cccc}
$(M_{21},\cdot)$&0&1&2&3\\
\hline
0&0&0&0&0\\
1&0&1&2&3\\
2&0&2&0&0\\
3&0&3&0&3\\
\multicolumn{5}{l}{\rule{0em}{1.5em}mult. $\cong M_{10}$}
\end{tabular}\quad
\begin{tabular}{r|cccc}
$(M_{22},\cdot)$&0&1&2&3\\
\hline
0&0&0&0&0\\
1&0&1&2&3\\
2&0&2&2&3\\
3&0&3&3&3\\
\multicolumn{5}{l}{\rule{0em}{1.5em}mult. $\cong M_{15}$}
\end{tabular}\quad
\begin{tabular}{r|cccc}
$(M_{23},\cdot)$&0&1&2&3\\
\hline
0&0&0&0&0\\
1&0&1&2&3\\
2&0&2&2&0\\
3&0&3&0&3\\
\multicolumn{5}{l}{\rule{0em}{1.5em}mult. $\cong M_{11}$}
\end{tabular}\\[1em]

\begin{tabular}{r|cccc}
$(M_{24},\cdot)$&0&1&2&3\\
\hline
0&0&0&0&0\\
1&0&1&2&3\\
2&0&2&1&3\\
3&0&3&3&3\\
\multicolumn{5}{l}{\rule{0em}{1.5em}mult. $\cong M_{16}$}
\end{tabular}\quad
\begin{tabular}{r|cccc}
$(M_{25},\cdot)$&0&1&2&3\\
\hline
0&0&0&0&0\\
1&0&1&2&3\\
2&0&2&2&0\\
3&0&3&0&3\\
\multicolumn{5}{l}{\rule{0em}{1.5em}mult. $\cong M_{11}$}
\end{tabular}\quad
\begin{tabular}{r|cccc}
$(M_{25},\cdot)$&0&1&2&3\\
\hline
0&0&0&0&0\\
1&0&1&2&3\\
2&0&2&1&3\\
3&0&3&3&0\\
\multicolumn{5}{l}{\rule{0em}{1.5em}mult. $\cong M_{12}$}
\end{tabular}\\[1em]

\begin{tabular}{r|cccc}
$(M_{25},\cdot)$&0&1&2&3\\
\hline
0&0&0&0&0\\
1&0&1&2&3\\
2&0&2&3&1\\
3&0&3&1&2\\
\multicolumn{5}{l}{\rule{0em}{1.5em}mult. $\cong M_{18}$}
\end{tabular}\quad
\begin{tabular}{r|cccc}
$(M_{26},\cdot)$&0&1&2&3\\
\hline
0&0&0&0&0\\
1&0&1&2&3\\
2&0&2&0&2\\
3&0&3&2&1\\
\multicolumn{5}{l}{\rule{0em}{1.5em}mult. $\cong M_{12}$}
\end{tabular}

\end{center}

\end{document}